\newif\ifpdf
\ifx\pdfoutput\undefined
\pdffalse 
\else
\pdfoutput=1 
\pdftrue \fi

\newif\ifarxiv
\arxivtrue

\documentclass[probth,final]{svjour}
\usepackage{cite}
\usepackage{amsmath}
\usepackage{amsfonts}
\usepackage{amssymb}
\usepackage{epsf}
\usepackage{graphicx}
\usepackage{graphics}
\usepackage{verbatim}
\usepackage{epsfig}

\IfFileExists{myowntimes.sty}
	{\usepackage{myowntimes}}
	{\usepackage{times}\usepackage{mathrsfs}}

\ifarxiv
\fi

\newenvironment{proofsect}[1]{\vskip0.1cm\noindent{\rmfamily\itshape #1.}}{\qed\vspace{0.15cm}}

\spnewtheorem{mylemma}[theorem]{Lemma}{\bf}{\it}
\spnewtheorem{myproposition}[theorem]{Proposition}{\bf}{\it}
\spnewtheorem{mycorollary}[theorem]{Corollary}{\bf}{\it}
\spnewtheorem{mydefinition}[theorem]{Definition}{\bf}{\it}
\spnewtheorem{myquestion}{Question}{\bf}{\it}
\spnewtheorem{myconjecture}[myquestion]{Conjecture}{\bf}{\it}
\numberwithin{equation}{section}
\numberwithin{theorem}{section}

\newcommand{\textd}{\text{\rm d}\mkern0.5mu}

\newcommand{\texte}{\text{\rm e}}

\newcommand{\1}{\operatorname{\sf 1}}

\newcommand{\DD}{\mathcal D}

\newcommand{\NN}{\mathcal N}

\newcommand{\prob}{\mathbb P}

\newcommand{\XX}{\mathcal X}

\newcommand{\B}{\mathbb B}

\newcommand{\E}{\mathbb E}

\newcommand{\G}{\mathbb G}
\newcommand{\BbbH}{\mathbb H}

\newcommand{\N}{\mathbb N}

\newcommand{\BbbP}{\mathbb P}

\newcommand{\R}{\mathbb R}

\newcommand{\Z}{\mathbb Z}

\newcommand{\scrB}{\mathscr{B}}
\newcommand{\scrC}{\mathscr{C}}

\newcommand{\scrF}{\mathscr{F}}
\newcommand{\scrG}{\mathscr{G}}

\newcommand{\scrS}{\mathscr{S}}

\newcommand{\scrW}{\mathscr{W} }
\newcommand{\scrX}{\mathscr{X}}

\newcommand{\twoeqref}[2]{(\ref{#1}--\ref{#2})}
\newcommand{\cc}{{\text{\rm c}}}

\def\c_#1{c_{\hbox{$\scriptscriptstyle #1$}}}

\def\myffrac#1#2 in #3{\raise 2.6pt\hbox{$#3 #1$}\mkern-1.5mu\raise 0.8pt\hbox{$#3/$}\mkern-1.1mu\lower 1.5pt\hbox{$#3 #2$}}
\newcommand{\ffrac}[2]{\mathchoice%
	{\myffrac{#1}{#2} in \scriptstyle}
	{\myffrac{#1}{#2} in \scriptstyle}
	{\myffrac{#1}{#2} in \scriptscriptstyle}
	{\myffrac{#1}{#2} in \scriptscriptstyle}
}

\begin{document}

\title{Quenched invariance principle for simple\\random walk on percolation clusters}
\titlerunning{Simple random walk on percolation clusters}
\author{Noam Berger\inst{1,2}
	\and Marek Biskup\inst{2}}
\authorrunning{Noam Berger and Marek Biskup}

\institute{Department of Mathematics, Caltech, Pasadena, CA 91125, U.S.A.
\and 
Department of Mathematics, UCLA, Los Angeles, CA 90095, U.S.A.}

\date{}
\maketitle

\renewcommand{\thefootnote}{}
\footnotetext{\vglue-0.41cm\footnotesize\copyright\,2005 by N.~Berger and M.~Biskup. Reproduction, by any means, of the entire article for non-commercial purposes is permitted without charge.}
\renewcommand{\thefootnote}{\arabic{footnote}}

\begin{abstract}
We consider the simple random walk on the (unique) infinite cluster of super-critical bond
percolation in~$\Z^d$ with~$d\ge2$. We prove that, for almost every percolation configuration, the path distribution of the walk converges weakly to that of non-degenerate, isotropic Brownian motion. Our analysis is based on the consideration of a harmonic deformation of the infinite cluster on which the random walk becomes a square-integrable martingale. The size of the deformation, expressed by the so called corrector, is estimated by means of ergodicity arguments.
\end{abstract}

\section{Introduction}
\subsection{Motivation and model}
\label{sec1.1}\noindent
Consider supercritical bond-percolation on $\Z^d$,~$d\ge2$, and the simple random walk on the (unique) infinite cluster. In~\cite{Sidoravicius-Sznitman} Sidoravicius and Sznitman asked the following question: Is it true
that for a.e.~configuration in which the origin belongs to the infinite
cluster, the random walk started at the origin exits the infinite symmetric slab $\{(x_1,\dots,x_d)\colon |x_d|\le N\}$
through the ``top'' side with probability tending to~$\ffrac 12$ as $N\to\infty$?
Sidoravicius and Sznitman managed to answer their question affirmatively in dimensions~$d\ge4$ but dimensions~$d=2,3$ remained open. In this paper we extend the desired conclusion to all~$d\ge2$.
As in~\cite{Sidoravicius-Sznitman}, we will do so by proving a quenched invariance principle for the paths of the walk.

Random walk on percolation clusters is only one of many instances of ``statistical mechanics in random media'' that have been recently considered by physicists and mathematicians.
Other pertinent examples include, e.g., various diluted spin systems, random copolymers~\cite{Soteros-Whittington}, spin glasses~\cite{Bolthausen-Sznitman,Talagrand}, random-graph models~\cite{Bollobas}, etc.
From this general perspective, the present problem is interesting for at least two reasons:
First, a good handle on simple random walk on a given graph is often a prerequisite for the understanding of more complicated processes, e.g., self-avoiding walk or loop-erased random walk. Second, information about the scaling properties of simple random walk on percolation cluster can, in principle, reveal some new important facts about the structure of the infinite cluster and/or its harmonic properties.

\smallskip
Let us begin developing the mathematical layout of the problem. Let~$\Z^d$ be the $d$-dimensional hypercubic lattice and let~$\B_d$ be the set of nearest neighbor edges. We will use~$b$ to denote a generic edge,~$\langle x,y\rangle$ to denote the edge between~$x$ and~$y$, and~$e$ to denote the edges from the origin to its nearest neighbors.
Let~$\Omega=\{0,1\}^{\B^d}$ be the space of all percolation configurations~$\omega=(\omega_b)_{b\in\B^d}$. Here~$\omega_b=1$ indicates that the edge~$b$ is occupied and~$\omega_b=0$ implies that it is vacant. Let~$\scrB$ be the Borel $\sigma$-algebra on~$\Omega$---defined using the product topology---and let~$\BbbP$ be an i.i.d.~measure such that~$\BbbP(\omega_b=1)=p$ for all~$b\in\B_d$. If~$x\overset\omega\longleftrightarrow\infty$ denotes the event 
that the site~$x$ 
belongs to an infinite self-avoiding 
path using only occupied bonds in~$\omega$, we write $\scrC_\infty=\scrC_\infty(\omega)$ for the set
\begin{equation}
\scrC_\infty(\omega)=\bigl\{x\in\Z^d\colon x\overset\omega\longleftrightarrow\infty\bigr\}.
\end{equation}
By Burton-Keane's uniqueness theorem~\cite{Burton-Keane}, the infinite cluster is unique and so $\scrC_\infty$ is connected with~$\BbbP$-probability one.

For each~$x\in\Z^d$, let~$\tau_x\colon\Omega\to\Omega$ be the ``shift by~$x$'' defined by~$(\tau_x\omega)_b=\omega_{x+b}$. Note that~$\BbbP$ is~$\tau_x$-invariant for all~$x\in\Z^d$.
Let~$p_\cc=p_\cc(d)$ denote the percolation threshold on~$\Z^d$ defined as the infimum of all~$p$'s for which $\BbbP(0\in\scrC_\infty)>0$. Let $\Omega_0=\{0\in\scrC_\infty\}$ and, for~$p>p_\cc$, define the measure~$\BbbP_0$ by
\begin{equation}
\BbbP_0(A)=\BbbP(A|\Omega_0),\qquad A\in\scrB.
\end{equation}
We will use~$\E_0$ to denote expectation with respect to~$\BbbP_0$.

For each configuration $\omega\in\Omega_0$, let $(X_n)_{n\ge0}$ be the simple random walk on $\scrC_\infty(\omega)$ started at the origin. Explicitly, $(X_n)_{n\ge0}$ is a Markov chain with state space~$\Z^d$, whose distribution~$P_{0,\omega}$ is defined by the transition probabilities
\begin{equation}
\label{E:1.3}
P_{0,\omega}(X_{n+1}=x+e|X_n=x)=\frac1{2d}\1_{\{\omega_e=1\}}\circ\tau_x,\qquad |e|=1,
\end{equation}
and
\begin{equation}
\label{E:1.3b}
P_{0,\omega}(X_{n+1}=x|X_n=x)=\sum_{e\colon |e|=1}\frac1{2d}\1_{\{\omega_e=0\}}\circ\tau_x,
\end{equation}
with the initial condition
\begin{equation}
\label{E:1.4}
P_{0,\omega}(X_0=0)=1.
\end{equation}
Thus, at each unit of time, the walk picks a neighbor at random and if the corresponding edge is occupied, the walk moves to this neighbor. If the edge is vacant, the move is suppressed.

\subsection{Main results}
Our main result is that for $\BbbP_0$-almost every $\omega\in\Omega_0$, the linear interpolation
of $(X_n)$, properly scaled, converges weakly to Brownian motion. For every~$T>0$, let~$(C[0,T],\scrW_T)$ be the space of continuous functions~$f\colon[0,T]\to\R$ equipped with the $\sigma$-algebra~$\scrW_T$ of Borel sets relative to the supremum topology. The precise statement is now as follows:

\begin{theorem}
\label{thm:mainthm}
Let~$d\ge2$, $p>p_\cc(d)$ and let~$\omega\in\Omega_0$. Let~$(X_n)_{n\ge0}$ be the random walk with law~$P_{0,\omega}$ and let
\begin{equation}
\label{1.5a}
\widetilde B_n(t)=\frac1{\sqrt n}\bigl(
X_{\lfloor tn\rfloor}+(tn-\lfloor tn\rfloor)(X_{\lfloor tn\rfloor+1}-X_{\lfloor tn\rfloor})\bigr),
\qquad t\ge0.
\end{equation}
Then for all~$T>0$ and for~$\BbbP_0$-almost every~$\omega$, the law of~$(\widetilde B_n(t)\colon 0\le t\le T)$ on $(C[0,T],\scrW_T)$ converges weakly to the law of an isotropic Brownian motion $(B_t\colon 0\le t\le T)$ whose diffusion constant, $D=E(|B_1|^2)>0$, depends only on the percolation parameter~$p$ and the dimension~$d$.
\end{theorem}

The Markov chain~$(X_n)_{n\ge0}$ represents only one of two natural ways to define a simple random walk on the supercritical percolation cluster. 
Another possibility is that, at each unit of time, the walk moves to a site chosen uniformly at random from the \emph{accessible} neighbors, i.e., the walk takes no pauses. In order to define this process, let~$(T_k)_{k\ge0}$ be the sequence of stopping times that mark the moments when the walk~$(X_n)_{n\ge0}$ made a move. Explicitly, $T_0=0$ and
\begin{equation}
\label{1.6a}
T_{n+1}=\inf\{k>T_n\colon X_k\ne X_{k-1}\}, \qquad n\ge0.
\end{equation}
Using these stopping times---which are~$P_{0,\omega}$-almost surely finite for all~$\omega\in\Omega_0$---we define a new Markov chain~$(X_n')_{n\ge0}$ by
\begin{equation}
\label{1.7a}
X_n'=X_{T_n},\qquad n\ge0.
\end{equation}
It is easy to see that~$(X_n')_{n\ge0}$ has the desired distribution. Indeed, the walk starts at the origin and its transition probabilities are given by
\begin{equation}
P_{0,\omega}(X_n'=x+e|X_n'=x)=\frac{\1_{\{\omega_e=1\}}\circ\tau_x}
{\sum_{e'\colon|e'|=1}\1_{\{\omega_{e'}=1\}}\circ\tau_x},\qquad |e|=1.
\end{equation}
A simple modification of the arguments leading to Theorem~\ref{thm:mainthm} allows us to establish a functional central limit theorem for this random walk as well:

\begin{theorem}
\label{thm:2ndmainthm}
Let~$d\ge2$, $p>p_\cc(d)$ and let~$\omega\in\Omega_0$. Let~$(X'_n)_{n\ge0}$ be the random walk defined from~$(X_n)_{n\ge0}$ as described in~\eqref{1.7a} and let $\widetilde B_n'(t)$ be the linear interpolation of~$(X_k')_{0\le k\le n}$ defined by~\eqref{1.5a} with~$(X_k)$ replaced by~$(X_k')$. Then for all~$T>0$ and for~$\BbbP_0$-almost every~$\omega$, the law of~$(\widetilde B_n'(t)\colon 0\le t\le T)$ on $(C[0,T],\scrW_T)$ converges weakly to the law of an isotropic Brownian motion $(B_t\colon 0\le t\le T)$ whose diffusion constant, $D'=E(|B_1|^2)>0$, depends only on the percolation parameter~$p$ and the dimension~$d$.
\end{theorem}

De Gennes~\cite{deGennes}, who introduced the problem of random walk on percolation cluster to the physics community, thinks of the walk as the motion of ``an ant in a labyrinth.'' From this perspective, the ``lazy'' walk~$(X_n)$ corresponds to a ``blind'' ant, while the ``agile'' walk $(X_n')$ represents a ``myopic'' ant. While the character of the scaling limit of the two ``ants'' is the same, there seems to be some distinction in the rate the scaling limit is approached, cf~\cite{Harris-at-al} and references therein.
As we will see in the proof, the diffusion constants~$D$ and~$D'$ are related via~$D'=D\Theta^2$, where~$\Theta^{-1}$ is the expected degree of the origin normalized by~$2d$, cf~\eqref{E:Theta}.

There is actually yet another way how to ``put'' simple random walk on~$\scrC_\infty$, and that is to use continuous time. Here the corresponding result follows by combining the CLT for the ``lazy'' walk with an appropriate Renewal Theorem for exponential waiting times.

\subsection{Discussion and related work}
The subject of random walk in random environment has a long history; we refer to, e.g.,~\cite{Bolthausen-Sznitman,Zeitouni} for recent overviews of (certain parts of) this field.
On general grounds, each random-media problem comes in two distinct flavors: \emph{quenched}, corresponding to the situations discussed above where the walk is distributed according to an~$\omega$-dependent measure~$P_{0,\omega}$, and \emph{annealed}, in which the path distribution of the walk is taken from the averaged measure~$A\mapsto\E_0(P_{0,\omega}(A))$.
Under suitable ergodicity assumptions, the annealed problem typically corresponds to the quenched problem averaged over the starting point. Yet the distinction is clear: In the annealed setting the slab-exit problem from Sect.~\ref{sec1.1} is trivial by the symmetries of the averaged measure, while its answer is \emph{a priori} very environment-sensitive in the quenched measure.

An annealed version of our theorems was proved in the 1980s by De Masi~\emph{et al}~\cite{demas1,demas2}, based on earlier results of Kozlov~\cite{Kozlov}, Kipnis and Varadhan~\cite{Kipnis-Varadhan} and others in the context of random walk in a field of random conductances.
(The results of~\cite{demas1,demas2} were primarily two-dimensional but, with the help of~\cite{Barlow}, they apply to all~$d\ge2$; cf~\cite{Sidoravicius-Sznitman}.)
A number of proofs of quenched invariance principles have appeared in recent years
for the cases where an annealed principle was already known. The most relevant paper is that of Sidoravicius and Sznitman~\cite{Sidoravicius-Sznitman} which established Theorem~\ref{thm:2ndmainthm} for random walk among random conductances in all~$d\ge1$ and, using a very different method, also for random walk on percolation in~$d\ge4$. 
(Thus our main theorem is new only in~$d=2,3$.) The~$d\ge4$ proof is based on the fact that two independent random walk paths will intersect only very little---something hard to generalize to~$d=2,3$. As this paper shows, the argument for random conductances is somewhat more flexible.

\begin{figure}[t]
{\epsfig{figure=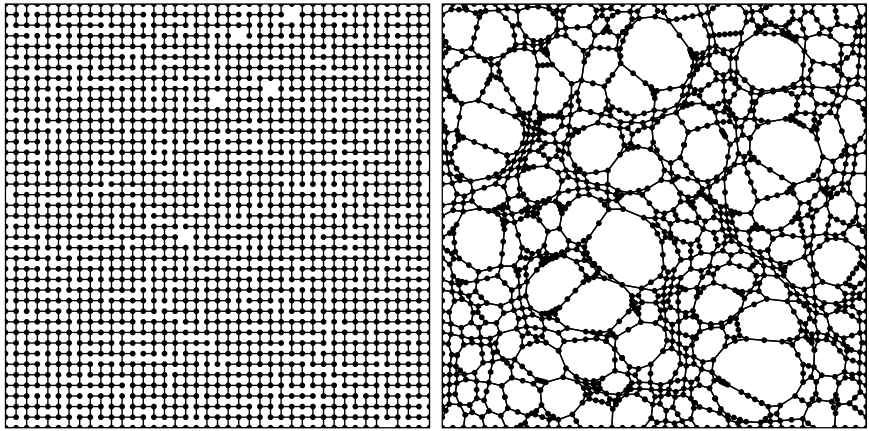, width=\textwidth}}
\smallskip
\caption{A portion of the infinite cluster~$\scrC_\infty=\scrC_\infty(\omega)$ before (left) and after (right) the harmonic deformation $x\mapsto x+\chi(x,\omega)$. Here~$p=0.75$ is already so large that all but a few sites in the entire block belong to~$\scrC_\infty$. Upon the deformation, all ``holes'' (i.e., dual connected components) get considerably stretched and rounded while the ``dangling ends'' collapse onto the rest of the structure.
}
\label{fig1}
\end{figure}  

Another paper of relevance is that of Rassoul-Agha and Sepp\"al\"ainen~\cite{rsagha} where a quenched invariance principle was established for \emph{directed} random walks in (space-time) random environments.
The directed setting offers the possibility to use independence more efficiently---every time step the walk enters a new environ\-ment---but the price to pay for this is the lack of reversibility.
The directed nature of the environment also permits consideration of distributions with a drift for which a CLT is not even expected to generally hold in the undirected setting; see~\cite{Berger-Gantert-Peres,Sznitman} for an example of ``pathologies'' that may arise.

Finally, there have been been a number of results dealing with harmonic properties of the simple random walk on percolation clusters. Grimmett, Kesten and Zhang~\cite{Grimmett-Kesten-Zhang} proved via ``electrostatic techniques'' that this random walk is transient in~$d\ge3$; extensions concerning the existence of various ``energy flows'' appeared in~\cite{Hoffman,Levin-Peres,Hoffman-Mossel,Mossel-Haggstrom,Angel-Benjamini-Berger-Peres}.
A great amount of effort has been spent on deriving estimates on the heat-kernel---i.e., the probability that the walk is at a particular site after~$n$ steps. 
The first such bounds were obtained by Heicklen and Hoffman~\cite{Heicklen-Hoffman}.
Later Mathieu and Remy~\cite{Mathieu-Remy} realized that the right way to approach heat-kernel estimates was through harmonic function theory of the infinite cluster and thus significantly improved the results of~\cite{Heicklen-Hoffman}. Finally, Barlow~\cite{Barlow} obtained, using again harmonic function theory, Gaussian upper and lower bounds for the heat kernel. We refer to~\cite{Barlow} for further references concerning this area of research.

\smallskip
\noindent\emph{Note}:
At the time a preprint version of this paper was first circulated, we learned that Mathieu and Piatnitski had announced a proof of the same result (albeit in continuous-time setting). Their proof, which has in the meantime been posted~\cite{Mathieu-Piatnitski}, is close in spirit to that of Theorem~1.1 of~\cite{Sidoravicius-Sznitman}; the main tools are Poincar\'e inequalities, heat-kernel estimates and homogenization theory.

\subsection{Outline}
Let us outline the main steps of our proof of Theorems~\ref{thm:mainthm} and~\ref{thm:2ndmainthm}. The principal idea---which permeates in various disguises throughout the work of Papanicolau and Varadhan~\cite{Papanicolau-Varadhan}, Kozlov~\cite{Kozlov}, Kipnis and Varadhan~\cite{Kipnis-Varadhan}, De Masi \emph{et al}~\cite{demas1,demas2}, Sidoravicius and Sznitman~\cite{Sidoravicius-Sznitman} and others---is to consider an embedding of~$\scrC_\infty(\omega)$ into the Euclidean space that makes the corresponding simple random walk a martingale. Formally, this is achieved by finding an~$\R^d$-valued discrete harmonic function on~$\scrC_\infty$ with a linear growth at infinity. The distance between the natural position of a site~$x\in\scrC_\infty$ and its counterpart in this \emph{harmonic embedding} is expressed in terms of the so-called \emph{corrector}~$\chi(x,\omega)$ which is a principal object of study in this paper. See Fig.~\ref{fig1} for an illustration. 

It is clear that the corrector can be defined in any finite volume by solving an appropriate discrete Dirichlet problem (this is how Fig.~1 was drawn); the difficult part is to define the corrector in infinite volume while maintaining the natural (distributional) invariance with respect to shifts of the underlying lattice. Actually, there is an alternative, probabilistic definition of the corrector,
\begin{equation}
\chi(x,\omega)=\lim_{n\to\infty}\bigl(E_{x,\omega}(X_n)-E_{0,\omega}(X_n)\bigr).
\end{equation}
However, the only proof we presently have for the existence of such a limit is by following, rather closely, the constructions from Sect.~\ref{sec2.3}.

Once we have the corrector under control, the proof splits into two parts: (1)~proving that the martingale---i.e., the walk on the deformed graph---converges to Brownian motion and (2)~proving that the deformation of the path caused by the change of embedding is negligible. The latter part (which is the principal contribution of this work) amounts to a sublinear bound on the corrector~$\chi(x,\omega)$ as a function of~$x$. Here, somewhat unexpectedly, our level of control is considerably better in~$d=2$ than in~$d\ge3$. In particular, our proof in~$d=2$ avoids using any of the recent sophisticated discrete-harmonic analyses but, to handle all~$d\ge2$ uniformly, we need to invoke the main result of Barlow~\cite{Barlow}.
The proof is actually carried out along these lines only for the setting in Theorem~\ref{thm:mainthm}; Theorem~\ref{thm:2ndmainthm} follows by noting that the time scales of both walks are comparable.

\smallskip
Here is a summary of the rest of this paper: 
In Sect.~\ref{sec:defcor} we introduce the aforementioned corrector and prove some of its basic properties.
Sect.~\ref{sec:erg} collects the needed facts about ergodic properties of the Markov chain ``on environments.'' 
Both sections are based on previously known material; proofs have been included to make the paper self-contained.
The novel parts of the proof---sublinear bounds on the corrector---appear in Sects.~\ref{sec:sec4}-\ref{sec:sublin}. The actual proofs of our main theorems are carried out in Sect.~\ref{sec:invprnp}.
The Appendix (Sects.~\ref{appA} and~\ref{appB}) contains the proof of an upper bound for the transition probabilities of our random walk, further discussion and some conjectures.

\begin{figure}[t]
\centerline{\epsfig{figure=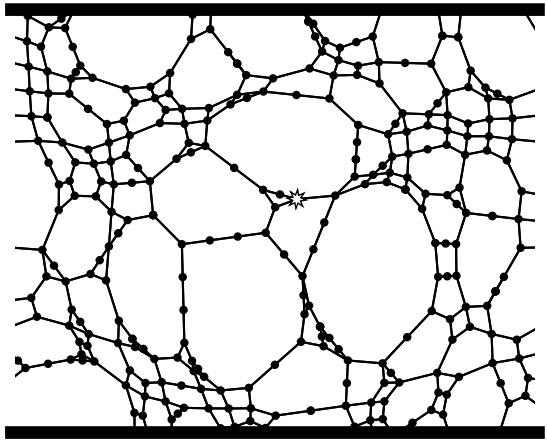, width=0.6\textwidth}}
\smallskip
\caption{The harmonic deformation of a percolation configuration in the symmetric slab $\{(x_1,x_2)\in\Z^2\colon |x_2|\le N\}$. The star denotes the new position of the origin which in the undeformed configuration was right in the center. The relative vertical shift of the origin corresponds to the deviation of~$P_{0,\omega}(\text{top hit before bottom})$ from one half. The figure also has an interesting electrostatic interpretation: If the bottom and top bars are set to potentials $-1$ and~$+1$, respectively, then the site with \emph{deformed} coordinates~$(x_1,x_2)$ has potential~$x_2/N$.
}
\label{fig2}
\end{figure}

\section{Corrector---construction and harmonicity}
\label{sec:defcor}\noindent
In this section we will define and study the aforementioned corrector which is the basic instrument of our proofs. The main idea is to consider the Markov chain ``on environments'' (Sect.~\ref{sec2.1}). The relevant properties of the corrector are listed in Theorem~\ref{thm2.1} (Sect.~\ref{sec2.2}); the proofs are based on spectral calculus (see Sect.~\ref{sec2.3}).

\subsection{Markov chain ``on environments''}
\label{sec2.1}\noindent
As is well known, cf~Kipnis and Varadhan~\cite{Kipnis-Varadhan}, the Markov chain~$(X_n)_{n\ge0}$ in \twoeqref{E:1.3}{E:1.4} induces a Markov chain on~$\Omega_0$, which can be interpreted as the trajectory of ``environments viewed from the perspective of the 
walk.'' The transition probabilities of this chain are given by the kernel $Q\colon\Omega_0\times\scrB\to[0,1]$,
\begin{equation}
\label{1.5}
Q(\omega,A)=\frac1{2d}\sum_{e\colon|e|=1}\bigl(\1_{\{\omega_{e}=1\}}\1_{\{\tau_{e}\omega\in A\}}
+\1_{\{\omega_{e}=0\}}\1_{\{\omega\in A\}}\bigr).
\end{equation}
Our basic observations about the induced Markov chain are as follows:

\begin{mylemma}
\label{lemma1.1}
For every bounded measurable~$f\colon\Omega\to\R$ and every~$e$ with~$|e|=1$,
\begin{equation}
\label{1.7}
\E_0\bigl(f\circ\tau_{e}\1_{\{\omega_{e}=1\}}\bigr)=\E_0\bigl(f\1_{\{\omega_{-e}=1\}}\bigr),
\end{equation}
where~$-e$ is the bond that is opposite to~$e$.
As a consequence, $\BbbP_0$ is reversible and, in particular, stationary for Markov kernel~$Q$.
\end{mylemma}

\begin{proofsect}{Proof}
First we will prove~\eqref{1.7}. Neglecting the normalization by~$\BbbP(0\in\scrC_\infty)$, we need that
\begin{equation}
\label{1.8}
\E\bigl(f\circ\tau_{e}\1_{\{0\in\scrC_\infty\}}\1_{\{\omega_{e}=1\}}\bigr)=\E\bigl(f\1_{\{0\in\scrC_\infty\}}\1_{\{\omega_{-e}=1\}}\bigr).
\end{equation}
This will follow from $\1_{\{\omega_{e}=1\}}=\1_{\{\omega_{-e}=1\}}\circ\tau_{e}$ and the fact that, on~$\{\omega_{e}=1\}$ we have $\1_{\{0\in\scrC_\infty\}}=\1_{\{0\in\scrC_\infty\}}\circ\tau_{e}$.
Indeed, these observations imply
\begin{equation}
f\circ\tau_{e}\1_{\{0\in\scrC_\infty\}}\1_{\{\omega_{e}=1\}}=\bigl(f\1_{\{0\in\scrC_\infty\}}\1_{\{\omega_{-e}=1\}}\bigr)\circ\tau_{e}
\end{equation}
and~\eqref{1.8} then follows by the shift invariance of~$\BbbP$.

From~\eqref{1.7} we deduce that for any bounded, measurable~$f,g\colon\Omega\to\R$,
\begin{equation}
\label{1.9}
\E_0\bigl(f(Qg)\bigr)=\E_0\bigl(g(Qf)\bigr),
\end{equation}
where~$Qf\colon\Omega\to\R$ is the function
\begin{equation}
\label{2.6a}
(Qf)(\omega)=\frac1{2d}\sum_{e\colon|e|=1}\bigl(\1_{\{\omega_{e}=1\}}f(\tau_{e}\omega)
+\1_{\{\omega_{e}=0\}}f(\omega)\bigr).
\end{equation}
Indeed, splitting the last sum into two terms, the second part reproduces exactly on both sides of~\eqref{1.9}. For the first part we apply~\eqref{1.7} and note that averaging over~$e$ allows us to neglect the negative sign in front of~$e$ on the right-hand side.
But~\eqref{1.9} is the \emph{definition} of reversibility and, setting~$f=1$ and noting that~$Q1=1$, we also get the stationarity of~$\BbbP_0$.
\end{proofsect}

Lemma~\ref{lemma1.1} underlines our main reason to work primarily with the ``lazy'' walk. For the ``agile'' walk, to get a stationary law on environments, one has to weigh~$\BbbP_0$ by the degree of the origin---a factor that would drag through the entire derivation.

\subsection{Kipnis-Varadhan construction}
\label{sec2.2}\noindent
Next we will adapt the construction of Kipnis and Varadhan~\cite{Kipnis-Varadhan} to the present situation.
Let~$L^2=L^2(\Omega,\scrB,\BbbP_0)$ be the space of all Borel-measurable, square integrable functions on~$\Omega$. Abusing the notation slightly, we will use ``$L^2$'' both for $\R$-valued functions as well as~$\R^d$-valued functions. We equip~$L^2$ with the inner product $(f,g)=\E_0(fg)$---with ``$fg$'' interpreted as the dot product of~$f$ and~$g$ when these functions are vector-valued. Let~$Q$ be the operator defined by~\eqref{2.6a}. Note that, when applied to a vector-valued function,~$Q$ acts like a scalar, i.e., independently on each component. 

From~\eqref{1.9} we know
\begin{equation}
(f,Qg)=(Qf,g)
\end{equation}
and so~$Q$ is symmetric.
An explicit calculation gives us
\begin{equation}
\begin{aligned}
\bigl|(f,Qf)\bigr|
&\le\frac1{2d}\sum_{e\colon|e|=1}\Bigl\{\bigl|(f,\1_{\{\omega_{e}=1\}} f\circ\tau_{e})\bigl|
+(f,\1_{\{\omega_{e}=0\}} f)\Bigr\}
\\
&\le\frac1{2d}\sum_{e\colon|e|=1}\Bigl\{\bigl|(f,\1_{\{\omega_{e}=1\}} f)\bigl|^{\ffrac12}
\bigl|(f,\1_{\{\omega_{-e}=1\}} f)\bigr|^{\ffrac12}
+(f,\1_{\{\omega_{e}=0\}} f)\Bigr\}
\\
&\le\frac1{2d}\sum_{e\colon|e|=1}\Bigl\{(f,\1_{\{\omega_{e}=1\}} f)
+(f,\1_{\{\omega_{e}=0\}} f)\Bigr\}=(f,f)
\end{aligned}
\end{equation}
and so~$\Vert Q\Vert_{L^2}\le1$. In particular,~$Q$ is self-adjoint and $\text{spec}(Q)\subset[-1,1]$.

Let~$V\colon\Omega\to\R^d$ be the local drift at the origin, i.e.,
\begin{equation}
V(\omega)=\frac1{2d}\sum_{e\colon|e|=1}e\1_{\{\omega_{e}=1\}}.
\end{equation}
(We will only be interested in~$V(\omega)$ for~$\omega\in\Omega_0$, but that is of no consequence here.) Clearly, since 
$V$ is bounded,
we have~$V\in L^2$. For each~$\epsilon>0$, let~$\psi_\epsilon\colon\Omega\to\R^d$ be the solution of
\begin{equation}
\label{2.4}
(1+\epsilon-Q)\psi_\epsilon=V.
\end{equation}
Since~$1-Q$ is a non-negative operator,~$\psi_\epsilon$ is well-defined and~$\psi_\epsilon\in L^2$ for all $\epsilon>0$.
The following theorem is the core of the whole theory:

\begin{theorem}
\label{thm2.1}
There is a function~$\chi\colon\Z^d\times\Omega_0\to\R^d$ such that for every~$x\in\Z^d$,
\begin{equation}
\label{2.5}
\lim_{\epsilon\downarrow0}\1_{\{x\in\scrC_\infty\}}(\psi_{\epsilon}\circ\tau_x-\psi_{\epsilon})\,
=\chi(x,\cdot),\qquad \text{\rm in }L^2.
\end{equation}
Moreover, the following properties hold:
\settowidth{\leftmargini}{(11)}
\begin{enumerate}
\item[(1)] (Shift invariance)
For $\BbbP_0$-almost every~$\omega\in\Omega_0$,
\begin{equation}
\label{2.12d}
\chi(x,\omega)-\chi(y,\omega)=\chi\bigl(x-y,\tau_y(\omega)\bigr)
\end{equation}
holds for all~$x,y\in\scrC_\infty(\omega)$.
\item[(2)] (Harmonicity)
For $\BbbP_0$-almost every~$\omega\in\Omega_0$, the function
\begin{equation}
x\mapsto \chi(x,\omega)+x
\end{equation}
is harmonic with respect to the transition probabilities~\twoeqref{E:1.3}{E:1.3b}.
\item[(3)] (Square integrability)
There exists a constant~$C<\infty$ such that
\begin{equation}
\label{eq:l2bnd}
\bigl\Vert[\chi(x+e,\cdot)-\chi(x,\cdot)]\1_{\{x\in\scrC_\infty\}}\1_{\{\omega_{e}=1\}}\circ\tau_x\bigr\Vert_2<C
\end{equation}
is true for all~$x\in\Z^d$ and all~$e$ with~$|e|=1$.
\end{enumerate}
\end{theorem}

The rest of this section is spent on proving Theorem~\ref{thm2.1}.
The proof is based on spectral calculus and it closely follows the corresponding arguments from~\cite{Kipnis-Varadhan}.
Alternative constructions invoke projection arguments, cf~\cite{Olla-notes,Mathieu-Piatnitski}.

\subsection{Spectral calculations}
\label{sec2.3}\noindent
Let~$\mu_V$ denote the spectral measure of~$Q\colon L^2\to L^2$ associated with function~$V$, i.e., for every bounded, continuous~$\Phi\colon[-1,1]\to\R$, we have
\begin{equation}
\label{2.15k}
\bigl(V,\Phi(Q)V\bigr)=\int_{-1}^1\Phi(\lambda)\mu_V(\textd \lambda).
\end{equation}
(Since~$Q$ acts as a scalar,~$\mu_V$ is the sum of the ``usual'' spectral measures for the Cartesian components of~$V$.)
In the integral we used that, since~$\text{spec}(Q)\in[-1,1]$, the measure~$\mu_V$ is supported entirely in~$[-1,1]$. The first observation, made already by Kipnis and Varadhan, is stated as follows:

\begin{mylemma}
\label{lemma2.2}
\begin{equation}
\int_{-1}^1\frac1{1-\lambda}\mu_V(\textd\lambda)<\infty.
\end{equation}
\end{mylemma}

\begin{proofsect}{Proof}
With some caution concerning the infinite cluster, the proof is a combination of arguments right before Theorem~1.3 of~\cite{Kipnis-Varadhan} and those in the proof of Theorem~4.1 of~\cite{Kipnis-Varadhan}. Let~$f\in L^2$ be a bounded real-valued function and note that, by Lemma~\ref{lemma1.1} and the symmetry of the sums,
\begin{equation}
\label{eq:2.17}
\sum_{e\colon|e|=1}e\,\E_0(f\1_{\{\omega_{e}=1\}})
=\frac12\sum_{e\colon|e|=1}e\,\E_0\bigl((f-f\circ\tau_{e})\1_{\{\omega_{e}=1\}}\bigr).
\end{equation}
Hence, for every~$a\in\R^d$ we get
\begin{equation}
\begin{aligned}
(f,a\cdot V)
&=\frac12\frac1{2d}
\sum_{e\colon|e|=1}(e\cdot a)\E_0\bigl((f-f\circ\tau_{e})\1_{\{\omega_{e}=1\}}\bigr)
\\
&\le\frac12\Bigl(\frac1{2d}\sum_{e\colon|e|=1}(e\cdot a)^2\BbbP(\omega_{e}=1)\Bigr)^{\ffrac12}
\\
&\qquad\qquad\quad\times
\Bigl(\frac1{2d}
\sum_{e\colon|e|=1}\E_0\bigl((f-f\circ\tau_{e})^2\1_{\{\omega_{e}=1\}}\bigr)\Bigr)^{\ffrac12}.
\end{aligned}
\end{equation}
The first term on the right-hand side equals a constant times~$|a|$, while Lemma~\ref{lemma1.1} allows us to rewrite the second term into
\begin{multline}
\label{2.21c}
\quad
\frac1{2d}
\sum_{e\colon|e|=1}\E_0\bigl((f-f\circ\tau_{e})^2\1_{\{\omega_{e}=1\}}\bigr)
\\
=2\frac1{2d}
\sum_{e\colon|e|=1}\E_0\bigl(f(f-f\circ\tau_{e})\1_{\{\omega_{e}=1\}}\bigr)=2\bigl(f,(1-Q)f\bigr).
\quad
\end{multline}
We thus get that there exists a constant~$C_0<\infty$ such that for all bounded~$f\in L^2$,
\begin{equation}
\label{eq:2.20}
\bigl|(f,a\cdot V)\bigr|\le C_0|a|\bigl(f,(1-Q)f\bigr)^{\ffrac12}.
\end{equation}
Applying~\eqref{eq:2.20} for~$f$ of the form~$f=a\cdot\Psi(Q)V$, summing~$a$ over coordinate vectors in~$\R^d$ and invoking \eqref{2.15k}, we find that for every bounded continuous $\Psi\colon[-1,1]\to\R$ and~$C=C_0\sqrt d$,
\begin{equation}
\left|\int\Psi(\lambda)\mu_V(\textd\lambda)\right|
\leq C\biggl(\,\int(1-\lambda)\Psi(\lambda)^2\mu_V(\textd\lambda)\biggr)^{\ffrac12}.
\end{equation}
Substituting~$\Psi_\epsilon(\lambda)=(\ffrac1\epsilon)\wedge\frac1{1-\lambda}$ for~$\Psi$ and noting that~$(1-\lambda)\Psi_\epsilon(\lambda)\le1$, we get
\begin{equation}
\int\Psi_\epsilon(\lambda)\mu_V(\textd\lambda)
\leq C\biggl(\,\int\Psi_\epsilon(\lambda)\mu_V(\textd\lambda)\biggr)^{\ffrac12}
\end{equation}
and so 
\begin{equation}
\int\Psi_\epsilon(\lambda)\mu_V(\textd\lambda)
\leq C^2.
\end{equation}
The Monotone Convergence Theorem now implies
\begin{equation}
\label{eq:intbdd}
\int\frac 1{1-\lambda}\mu_V(\textd\lambda)
=\sup_{\epsilon>0}\,\int\Psi_\epsilon(\lambda)\mu_V(\textd\lambda)
\le C^2<\infty,
\end{equation}
proving the desired claim.
\end{proofsect}

Using spectral calculus we will now prove:

\begin{mylemma}
\label{lemma2.3}
Let~$\psi_\epsilon$ be as defined in~\eqref{2.4}. Then
\begin{equation}
\label{2.14}
\lim_{\epsilon\downarrow0}\epsilon\Vert\psi_\epsilon\Vert_2^2=0.
\end{equation}
Moreover, for~$e$ with~$|e|=1$ let $G_{e}^{(\epsilon)}=\1_{\{0\in\scrC_\infty\}}\1_{\{\omega_{e}=1\}}(\psi_\epsilon\circ\tau_{e}-\psi_\epsilon)$. Then for all~$x\in\Z^d$ and all~$e$ with~$|e|=1$,
\begin{equation}
\label{2.15}
\lim_{\epsilon_1,\epsilon_2\downarrow0}\,\bigl\Vert G_{e}^{(\epsilon_1)}\circ\tau_x-G_{e}^{(\epsilon_2)}\circ\tau_x\bigr\Vert_2=0
\end{equation}
\end{mylemma}

\begin{proofsect}{Proof}
The main ideas are again taken more or less directly from the proof of Theorem~1.3 in~\cite{Kipnis-Varadhan}; some caution is necessary regarding the containment in the infinite cluster in the proof of~\eqref{2.15}.
By the definition of~$\psi_\epsilon$,
\begin{equation}
\epsilon\Vert\psi_\epsilon\Vert_2^2=\int_{-1}^1\frac\epsilon{(1+\epsilon-\lambda)^2}\mu_V(\textd\lambda).
\end{equation}
The integrand is dominated by~$\frac1{1-\lambda}$ and tends to zero as~$\epsilon\downarrow0$ for every~$\lambda$ in the support of~$\mu_V$. Then~\eqref{2.14} follows by the Dominated Convergence Theorem.

The second part of the claim is proved similarly: First we get rid of the~$x$-dependence by noting that, due to the fact that $G_{e}^{(\epsilon)}\circ\tau_x\ne0$ enforces~$x\in\scrC_\infty$, the translation invariance of~$\BbbP$ implies
\begin{equation}
\bigl\Vert G_{e}^{(\epsilon_1)}\circ\tau_x-G_{e}^{(\epsilon_2)}\circ\tau_x\bigr\Vert_2
\le
\bigl\Vert G_{e}^{(\epsilon_1)}-G_{e}^{(\epsilon_2)}\bigr\Vert_2.
\end{equation}
Next we square the right-hand side and average over all~$e$. Using that~$G_{e}\ne0$ also enforces~$\omega_{e}=1$ and applying~\eqref{eq:2.17}, we thus get
\begin{equation}
\frac1{2d}\sum_{e\colon|e|=1}\bigl\Vert G_{e}^{(\epsilon_1)}-G_{e}^{(\epsilon_2)}\bigr\Vert_2^2
=2\bigl(\psi_{\epsilon_1}-\psi_{\epsilon_2},(1-Q)(\psi_{\epsilon_1}-\psi_{\epsilon_2})\bigr).
\end{equation}
Now we calculate
\begin{multline}
\qquad
\bigl(\psi_{\epsilon_1}-\psi_{\epsilon_2},(1-Q)(\psi_{\epsilon_1}-\psi_{\epsilon_2})\bigr)
\\=
\int_{-1}^1\frac{(\epsilon_1-\epsilon_2)^2(1-\lambda)}{(1+\epsilon_1-\lambda)^2(1+\epsilon_2-\lambda)^2}
\mu_V(\textd\lambda).
\qquad
\end{multline}
The integrand is again bounded by~$\frac1{1-\lambda}$, for all~$\epsilon_1,\epsilon_2>0$, and it tends to zero as~$\epsilon_1,\epsilon_2\downarrow0$. The claim follows by the Dominated Convergence Theorem.
\end{proofsect}

Now we are ready to prove Theorem~\ref{thm2.1}:

\begin{proofsect}{Proof of Theorem~\ref{thm2.1}}
Let~$G_{e}^{(\epsilon)}\circ\tau_x$ be as in Lemma~\ref{lemma2.3}. Using~\eqref{2.15} we know that~$G_{e}^{(\epsilon)}\circ\tau_x$ converges in~$L^2$ as~$\epsilon\downarrow0$. We denote~$G_{x,x+e}=\lim_{\epsilon\downarrow0}G_{e}^{(\epsilon)}\circ\tau_x$.
Since $G_{e}^{(\epsilon)}\circ\tau_x$ is a gradient field on~$\scrC_\infty$, we have $G_{x,x+e}(\omega)+G_{x+e,x}(\omega)=0$ and, more generally,~$\sum_{k=0}^{n}G_{x_k,x_{k+1}}=0$ whenever~$(x_0,\dots,x_n)$ is a closed loop on~$\scrC_\infty$. Thus, we may define
\begin{equation}
\chi(x,\omega)\,\overset{\text{\rm def}}=\,\sum_{k=0}^{n-1}G_{x_k,x_{k+1}}(\omega),
\end{equation}
where $(x_0,x_1,\dots,x_n)$ is a nearest-neighbor path on~$\scrC_\infty(\omega)$ connecting~$x_0=0$ to~$x_n=x$. By the above ``loop'' conditions, the definition is independent of this path for almost every~$\omega\in\{x\in\scrC_\infty\}$.
The shift invariance~\eqref{2.12d} now follows from this definition and~$G_{x,x+e}=G_{0,e}\circ\tau_x$.

In light of shift invariance, to prove the harmonicity of~$x\mapsto x+\chi(x,\omega)$ it suffices to show that, almost surely,
\begin{equation}
\label{2.29a}
\frac1{2d}\sum_{e\colon|e|=1}\bigl[\chi(0,\cdot)-\chi(e,\cdot)\bigr]\1_{\{\omega_e=1\}}=V.
\end{equation}
Since~$\chi(e,\cdot)-\chi(0,\cdot)=G_{0,e}$, the left hand side is the $\epsilon\downarrow0$ limit of
\begin{equation}
\frac1{2d}\sum_{e\colon|e|=1}\bigl[\psi_{\epsilon}-\psi_{\epsilon}\circ\tau_e\bigr]\1_{\{\omega_e=1\}}=
(1-Q)\psi_{\epsilon}.
\end{equation}
The definition of~$\psi_\epsilon$ tells us that $(1-Q)\psi_{\epsilon}=-\epsilon\psi_{\epsilon}+V$.
From here we get~\eqref{2.29a} by recalling that~$\epsilon\psi_{\epsilon}(\omega)$ tends to zero in~$L^2$.

To prove the square integrability in part~(3) we note that, by the construction of the corrector,
\begin{equation}
\bigl[\chi(x+e,\cdot)-\chi(x,\cdot)\bigr]\1_{\{x\in\scrC_\infty\}}\1_{\{\omega_{e}=1\}}\circ\tau_x=G_{x,x+e}.
\end{equation}
But~$G_{x,x+e}$ is the $L^2$-limit of~$L^2$-functions~$G_{e}^{(\epsilon)}\circ\tau_x$ whose $L^2$-norm is bounded by that of~$G_{e}^{(\epsilon)}$. Hence~\eqref{eq:l2bnd} follows with~$C=\max_{e\colon|e|=1}\Vert G_{0,e}\Vert_2$.
\end{proofsect}

\section{Ergodic-theory input}
\label{sec:erg}
\noindent
Here we will establish some basic claims whose common feature is the use of ergodic theory. Modulo some care for the containment in the infinite cluster, all of these results are quite standard and their proofs (cf~Sect.~\ref{sec3.2}) may be skipped on a first reading. Readers interested only in the principal conclusions of this section should focus their attention on Theorems~\ref{thm3.1} and~\ref{thm3.2}.

\subsection{Statements}
\label{sec3.1}\noindent
Our first result concerns the convergence of ergodic averages for the Markov chain on environments. The claim that will suffice for our later needs is as follows:

\begin{theorem}
\label{thm3.1}
Let~$f\in L^1(\Omega,\scrB,\BbbP_0)$. Then for~$\BbbP_0$-almost all~$\omega\in\Omega$,
\begin{equation}
\lim_{n\to\infty}\,\frac1n\sum_{k=0}^{n-1}f\circ\tau_{X_k}(\omega)=\E_0(f),\qquad P_{0,\omega}\text{\rm-almost surely}.
\end{equation}
Similarly, if~$f\colon\Omega\times\Omega\to\R$ is measurable with $\E_0(E_{0,\omega}|f(\omega,\tau_{X_1}\omega)|)<\infty$, then
\begin{equation}
\lim_{n\to\infty}\,\frac1n\sum_{k=0}^{n-1}f(\tau_{X_k}\omega,\tau_{X_{k+1}}\omega)=\E_0\bigl(E_{0,\omega}(f(\omega,\tau_{X_1}\omega))\bigr)
\end{equation}
for~$\BbbP_0$-almost all~$\omega$ and~$P_{0,\omega}$-almost all trajectories of~$(X_k)_{k\ge0}$.
\end{theorem}

\smallskip
The next principal result of this section will be the ergodicity of the ``induced shift'' on~$\Omega_0$. To define this concept, let~$e$ be a vector with~$|e|=1$ and, for every $\omega\in\Omega_0$, let
\begin{equation}
\label{3.1}
n(\omega)=\min\bigl\{k>0\colon ke\in\scrC_\infty(\omega)\bigr\}.
\end{equation}
By Birkhoff's Ergodic Theorem we know that~$\{k>0\colon ke\in\scrC_\infty\}$ has positive density in~$\N$ and so $n(\omega)<\infty$ almost surely. Therefore we can define the map~$\sigma_e\colon\Omega_0\to\Omega_0$ by
\begin{equation}
\sigma_e(\omega)=\tau_{n(\omega)e}\,\omega.
\end{equation} 
We call~$\sigma_e$ the \emph{induced shift}. Then we claim:

\begin{theorem}
\label{thm3.2}
For every~$e$ with~$|e|=1$, the induced shift $\sigma_e\colon\Omega_0\to\Omega_0$ is $\BbbP_0$-preserving and ergodic with respect to~$\BbbP_0$.
\end{theorem}

Both theorems will follow once we establish of ergodicity of the Markov chain on environments (see Proposition~\ref{prop3.5}). For finite-state (irreducible) Markov chains the proof of ergodicity is a standard textbook material (cf~\cite[page~51]{Petersen}), but our state space is somewhat large and so alternative arguments are necessary. Since we could not find appropriate versions of all needed claims in the literature, we include complete proofs.

\subsection{Proofs}
\label{sec3.2}\noindent
We begin by Theorem~\ref{thm3.2} which will follow from a more general statement, Lemma~\ref{claim:genergd}, below.
Let $(\XX,\scrX,\mu)$ be a probability space, and let $T\colon\XX\to\XX$
be invertible, measure preserving and ergodic with respect to~$\mu$.
Let $A\in\scrX$ be of positive
measure, and define $n\colon A\to\N\cup\{\infty\}$ by
\begin{equation}
n(x)=\min\bigl\{k>0\colon T^k(x)\in A\bigr\}.
\end{equation}
The Poincar\'e Recurrence Theorem (cf~\cite[Sect.~2.3]{Petersen}) tells us that $n(x)<\infty$ almost surely. Therefore we can define, up to a set of measure zero, the 
map $S\colon A\to A$ by
\begin{equation}
S(x)=T^{n(x)}(x),\qquad x\in A.
\end{equation}
Then we have:

\begin{mylemma}
\label{claim:genergd}
$S$ is measure preserving and ergodic with respect to $\mu(\cdot|A)$.
It is also almost surely invertible with respect to the same measure.
\end{mylemma}

\begin{proofsect}{Proof}
(1)~$S$ is measure preserving:
For $j\ge1$, let $A_j=\{x\in A\colon n(x)=j\}$. Then the $A_j$'s are disjoint
and $\mu(A\setminus\bigcup_{j\geq 1}A_j)=0$. First we show that
\begin{equation}
\label{3.4}
i\neq j\quad\Rightarrow\quad
S(A_i)\cap S(A_j)=\emptyset.
\end{equation}
To do this, we use the fact that $T$ is invertible. Indeed, if~$x\in S(A_i)\cap S(A_j)$ for $1\le i<j$, then~$x=T^i(y)=T^j(z)$ for some~$y,z\in A$ with~$n(y)=i$ and~$n(z)=j$. But the fact that~$T$ is invertible implies that~$y=T^{j-i}(z)$, which means~$n(z)\le j-i<j$, a contradiction.
To see that~$S$ is measure preserving, we note that the restriction of~$S$ to~$A_j$ is~$T^j$, which is measure preserving. Hence,~$S$ is measure preserving on~$A_j$ and, by~\eqref{3.4}, on the disjoint union~$\bigcup_{j\geq 1}A_j$ as well.

(2)~$S$ is almost surely invertible: $S^{-1}(\{x\})\cap\{S\text{ is well defined}\}$ is a one-point set by the fact that~$T$ is itself invertible.

(3)~$S$ is ergodic:
Let $B\in\scrX$ be such that $B\subseteq A$ and $0<\mu(B)<\mu(A)$. Assume that~$B$ is
$S$-invariant. Then $S^{n}(x)\notin A\setminus B$ for all $x\in B$ and all $n\ge1$.
This means that for every $x\in B$ and every $k\ge1$ such that $T^{k}(x)\in A$,
we have $T^{k}(x)\notin A\setminus B$. If follows that $C=\bigcup_{k\ge1}T^k(B)$ is (almost-surely) $T$-invariant and $\mu(C)\in(0,1)$, a contradiction with the ergodicity of $T$.
\end{proofsect}

\begin{proofsect}{Proof of Theorem~\ref{thm3.2}}
We know that the shift $\tau_{e}$ is invertible, measure preserving and ergodic
with respect to~$\BbbP$. By Lemma~\ref{claim:genergd} the induced shift~$\sigma_e\colon\Omega_0\to\Omega_0$ is~$\BbbP_0$-preserving, almost-surely invertible and ergodic with respect to~$\BbbP_0$.
\end{proofsect}

In the present circumstances, Theorem~\ref{thm3.2} has one important consequence:

\begin{mylemma}
\label{lemma3.3}
Let~$B\in\scrB$ be a subset of~$\Omega_0$ such that for almost all~$\omega\in B$,
\begin{equation}
\label{3.5}
P_{0,\omega}(\tau_{X_1}\omega\in B)=1.
\end{equation}
Then~$B$ is a zero-one event under~$\BbbP_0$.
\end{mylemma}

\begin{proofsect}{Proof}
The Markov property and \eqref{3.5} imply that~$P_{0,\omega}(\tau_{X_n}\omega\in B)=1$ for all~$n\ge1$ and~$\BbbP_0$-almost every~$\omega\in B$. We claim that~$\sigma_e(\omega)\in B$ for~$\BbbP_0$-almost every~$\omega\in B$. Indeed, let~$\omega\in B$ be such that~$\tau_{X_n}\omega\in B$ for all~$n\ge1$, $P_{0,\omega}$-almost surely. Let~$n(\omega)$ be as in~\eqref{3.1} and note that we have~$n(\omega)e\in\scrC_\infty$. By the uniqueness of the infinite cluster, there is a path of finite length connecting~$0$ and~$n(\omega)e$. If~$\ell$ is the length of this path, we have~$P_{0,\omega}(X_\ell=n(\omega)e)>0$. This means that~$\sigma_e(\omega)=\tau_{n(\omega)e}(\omega)\in B$, i.e.,~$B$ is almost surely~$\sigma_e$-invariant. By the ergodicity of the induced shift,~$B$ is a zero-one event.
\end{proofsect}

Our next goal will be to prove that the Markov chain on environments is ergodic. Let $\XX=\Omega^\Z$ and define~$\scrX$ to be the product~$\sigma$-algebra on~$\XX$; $\scrX=\scrB^{\otimes\Z}$.
The space~$\XX$ is a space of two-sided sequences~$(\dots,\omega_{-1},\omega_0,\omega_1,\dots)$---the trajectories of the Markov chain on environments. (Note that the index on~$\omega$ is an index in the sequence which is unrelated to the value of the configuration at a point.) 
Let~$\mu$ be the measure on $(\XX,\scrX)$ such that for any~$B\in\scrB^{2n+1}$,
\begin{multline}
\qquad
\mu\bigl((\omega_{-n},\dots,\omega_n)\in B\bigr)
\\=\int_B \BbbP_0(\textd \omega_{-n})Q(\omega_{-n},\textd \omega_{-n+1})\cdots
Q(\omega_{n-1},\textd \omega_{n}),
\qquad
\end{multline}
where~$Q$ is the Markov kernel defined in~\eqref{1.5}. (Since~$\BbbP_0$ is preserved by~$Q$, these finite-dimensional measures are consistent and~$\mu$ exists and is unique by Kolmogorov's Theorem.) 
Clearly,~$(\tau_{X_k}(\omega))_{k\ge0}$ has the same law in~$\E_0(P_{0,\omega}(\cdot))$ as~$(\omega_0,\omega_1,\dots)$ has in~$\mu$.
Let $T\colon\XX\to\XX$ be the shift defined by~$(T\omega)_n=\omega_{n+1}$. Then $T$ is measure preserving.
 
\begin{myproposition}
\label{prop3.5}
$T$ is ergodic with respect to~$\mu$.
\end{myproposition}

\begin{proofsect}{Proof}
Let~$E_\mu$ denote expectation with respect to~$\mu$.
Pick $A\subseteq\XX$ that is measurable and $T$-invariant.
We need to show that 
\begin{equation}\label{eq:Aistriv}
\mu(A)\in\{0,1\}.
\end{equation}
Let $f\colon \Omega\to\R$ be defined as $f(\omega_0)=E_\mu(\1_A|\omega_0)$. First we claim that~$f=\1_A$ almost surely. Indeed, since $A$ is $T$-invariant, there exist $A_+\in\sigma(\omega_k\colon k>0)$ and $A_-\in\sigma(\omega_k\colon k<0)$ such that~$A$ and~$A_\pm$ differ only by null sets from one another. (This follows by approximation of~$A$ by finite-dimensional events and using the $T$-invariance of~$A$.) Now conditional on~$\omega_0$, the event~$A_+$ is independent of~$\sigma(\omega_k\colon k<0)$ and so L\'evy's Martingale Convergence Theorem gives us
\begin{equation}
\begin{aligned}
E_\mu(\1_A|\omega_0)&=E_\mu(\1_{A_+}|\omega_0)=E_\mu(\1_{A_+}|\omega_0,\omega_{-1},\ldots,\omega_{-n})
\\*[2mm]
&=E_\mu(\1_{A_-}|\omega_0,\omega_{-1},\ldots,\omega_{-n})\,\underset{n\to\infty}\longrightarrow\,
\1_{A_-}=\1_A,
\end{aligned}
\end{equation}
with all equalities valid~$\mu$-almost surely.

Next let $B\subset\Omega$ be defined by~$B=\{\omega_0\colon f(\omega_0)=1\}$. Clearly,~$B$ is~$\scrB$-measurable and, since the $\omega_0$-marginal of~$\mu$ is~$\BbbP_0$,
\begin{equation}
\mu(A)=E_\mu(f)=\BbbP_0(B).
\end{equation}
Hence, to prove~\eqref{eq:Aistriv}, we need to show that 
\begin{equation}
\label{eq:Bistriv}
\BbbP_0(B)\in\{0,1\}.
\end{equation}
But~$A$ is~$T$-invariant and so, up to sets of measure zero, if~$\omega_0\in B$ then~$\omega_1\in B$. This means that~$B$ satisfies condition~\eqref{3.5} of Lemma~\ref{lemma3.3} and so \eqref{eq:Bistriv} holds.
\end{proofsect}

Now we can finally prove Theorem~\ref{thm3.1}:

\begin{proofsect}{Proof of Theorem~\ref{thm3.1}}
Recall that~$(\tau_{X_k}(\omega))_{k\ge0}$ has the same law in~$\E_0(P_{0,\omega}(\cdot))$ as~$(\omega_0,\omega_1,\dots)$ has in~$\mu$. Hence, if~$g(\dots,\omega_{-1},\omega_0,\omega_1,\dots)=f(\omega_0)$ then
\begin{equation}
\lim_{n\to\infty}\frac1n\sum_{k=0}^\infty f\circ\tau_{X_k} \overset{\DD}= \lim_{n\to\infty}\frac1n\sum_{k=0}^\infty g\circ T^k.
\end{equation}
The latter limit exists by Birkhoff's Ergodic Theorem and (by Proposition~\ref{prop3.5}) equals $E_\mu(g)=\E_0(f)$ almost surely. The second part is proved analogously.
\end{proofsect}

\section{Sublinearity along coordinate directions}
\label{sec:sec4}
\noindent
Equipped with the tools from the previous two sections, we can start addressing the main problem of our proof: the sublinearity of the corrector. Here we will prove the corresponding claim along the coordinate directions in~$\Z^d$.

\smallskip
Fix~$e$ with~$|e|=1$ and let~$n(\omega)$ be as defined in~\eqref{3.1}. Define a sequence~$n_k(\omega)$ inductively by $n_1(\omega)=n(\omega)$ and $n_{k+1}(\omega)=n_k(\sigma_e(\omega))$. The numbers~$(n_k)$, which are well-defined and finite on a set of full~$\BbbP_0$-measure, represent the successive ``arrivals'' of~$\scrC_\infty$ to the positive part of the coordinate axis in direction~$e$. Let~$\chi$ be the corrector defined in Theorem~\ref{thm2.1}. The main goal of this section is to prove the following theorem:

\begin{theorem}
\label{thm4.1}
For $\BbbP_0$-almost all~$\omega\in\Omega_0$,
\begin{equation}
\lim_{k\to\infty}\frac{\chi(n_k(\omega)e,\omega)}{k}=0.
\end{equation}
\end{theorem}

The proof is based on the following facts about the moments of~$\chi(n_k(\omega)e,\omega)$:

\begin{myproposition}
\label{lemma:next_corrector}
Abbreviate $v_e=v_e(\omega)=n_1(\omega)e$.  Then
\settowidth{\leftmargini}{(11)}
\begin{enumerate}
\item[(1)]
$\E_0(|\chi(v_e,\cdot)|)<\infty$.
\item[(2)]
$\E_0(\chi(v_e,\cdot))=0$.
\end{enumerate}
\end{myproposition}

The proof of this proposition will in turn be based on a bound on the tails of the length of the shortest path connecting the origin to~$v_e$.
We begin by showing that~$|v_e|$ has exponential tails:

\begin{mylemma}
\label{lemma-tails}
For each~$p>p_\cc$ there exists a constant~$a=a(p)>0$ such that for all~$e$ with~$|e|=1$,
\begin{equation}
\label{4.2k}
\BbbP_0\bigl(|v_e|>n\bigr)\le\texte^{-an},\qquad n\ge1.
\end{equation}
\end{mylemma}

\begin{proofsect}{Proof}
The proof uses a different argument in~$d=2$ and~$d\ge3$. In~$d\ge3$, we will use the fact that the slab-percolation threshold coincides with~$p_\cc$, as was proved by Grimmett and Marstrand~\cite{Grimmett-Marstrand}. Indeed, given~$p>p_\cc$, let~$K\ge1$ be so large that~$\Z^{d-1}\times\{1,\dots,K\}$ contains an infinite cluster almost surely. By the uniqueness of the percolation cluster in~$\Z^d$, this slab-cluster is almost surely a subset of~$\scrC_\infty$. Our bound in \eqref{4.2k} is derived as follows: Let~$A_K$ be the event that at least one of the sites in~$\{je\colon j=1,\dots,K\}$ is contained in the infinite connected component in~$\Z^{d-1}\times\{1,\dots,K\}$. Then~$\{|v_e|\ge Kn\}\cap\{0\in\scrC_\infty\}\subset\bigcap_{\ell\le n}\tau_{\ell Ke}(A)$. Since the events~$\tau_{\ell Ke}(A)$, $\ell=1,\dots,n$, are independent, letting~$p_K=\BbbP(A_K)$ we have
\begin{equation}
\BbbP(|v_e|\ge Kn,\,0\in\scrC_\infty)\le(1-p_K)^n,\qquad n\ge1.
\end{equation}
From here \eqref{4.2k} follows by choosing~$a$ appropriately.

In dimension~$d=2$, we will instead use a duality argument. Let~$\Lambda_n$ be the box $\{1,\dots,n\}\times\{1,\dots,n\}$. On $\{|v_e|\ge n\}\cap\{0\in\scrC_\infty\}$, none of the boundary sites $\{je\colon j=1,\dots,n\}$ are in~$\scrC_\infty$. So either at least one of these sites is in a finite component of size larger than~$n$ or there exists a dual crossing of~$\Lambda_n$ in the direction of~$e$. By the exponential decay of truncated connectivities (Theorem~8.18 of Grimmett~\cite{Grimmett}) and dual connectivities (Theorem~6.75 of Grimmett~\cite{Grimmett}), the probability of each of these events decays exponentially with~$n$.
\end{proofsect}

Our next lemma provides the requisite tail bound for the length of the shortest path between the origin and~$v_e$:

\begin{mylemma}
\label{claim:pathlen}
Let $L=L(\omega)$ be the length of the shortest occupied path from~$0$ to~$v_e$. Then there exist a constant~$C<\infty$ and $a>0$ such that for every $n\ge1$,
\begin{equation}
\BbbP_0(L>n)<C\texte^{-an}.
\end{equation}
\end{mylemma}

\begin{proofsect}{Proof}
Let~$\textd_\omega(0,x)$ be the length of the shortest path from~$0$ to~$x$ in configuration~$\omega$. Pick~$\epsilon>0$ such that~$\epsilon n$ is an integer. Then
\begin{equation}
\label{4.3a}
\{L>n\}\subset\bigl\{|v_e|\ge \epsilon n\bigr\}\cup\bigcup_{k=1}^{\epsilon n}\bigl\{\textd_\omega(0,ke)>n;\,0,ke\in\scrC_\infty\bigr\}.
\end{equation}
In light of Lemma~\ref{lemma-tails}, the claim will follow once we show that the probability of all events in the giant union on the right-hand side is bounded by~$\texte^{-a'n}$ with some~$a'>0$ (independently of~$k$).

We will use the following large-deviation result from Theorem~1.1 of Antal and Pisztora~\cite{Antal-Pisztora}: There exist constants~$a,\rho<\infty$ such that
\begin{equation}
\label{4.2}
\BbbP\bigl(\textd_\omega(0,x)>\rho|x|\bigr)\le \texte^{-a|x|}
\end{equation}
once~$|x|$ is sufficiently large. Unfortunately, we cannot use this bound in~\eqref{4.3a} directly, because~$ke$ can be arbitrarily close to~$0$ (in $\ell^\infty$ distance on~$\Z^d$). To circumvent this problem, let~$w_e$ be the site~$-me$ such that~$m=\min\{m'>\epsilon n\colon -m'e\in\scrC_\infty\}$ and let~$A_{x,y}=\{\textd_\omega(x,y)\ge\ffrac n2,\,x,y\in\scrC_\infty\}$. Then, on~$\{\textd_\omega(0,x)>n\}$, either~$|w_e|>2\epsilon n$ or at least one site ``between''~$-2\epsilon n e$ and~$-\epsilon ne$ is connected to either~$0$ or~$ke$ by a path longer than~$\ffrac n2$. Since on~$\{|w_e|>2\epsilon n\}$ we must have~$|v_{-e}\circ\sigma_{-e}^m|>\epsilon n$ for at least one~$m=1,\dots\epsilon n$, we have
\begin{multline}
\quad
\bigl\{\textd_\omega(0,ke)>n;\,0,ke\in\scrC_\infty\bigr\}
\\
\subset\Bigl(\bigcup_{m=1}^{\epsilon n}\sigma_{-e}^m\bigl(\{|v_{-e}|\ge\epsilon n\}\bigr)
\cup\bigcup_{\epsilon n\le \ell\le2\epsilon n}\bigl(A_{0,-\ell e}\cup A_{ke,-\ell e}\bigr)\Bigr).
\quad
\end{multline}
Now all events in the first giant union have the same probability, which is exponentially small by Lemma~\ref{lemma-tails}. As to the second union, by~\eqref{4.2} we know that
\begin{equation}
\BbbP_0(A_{0,-\ell e})\le \texte^{-a\ell}\le \texte^{-a\epsilon n}
\end{equation}
whenever~$\epsilon$ is so small that~$4\epsilon\rho\le1$, and a similar bound holds for~$A_{ke,-\ell e}$ as well (except that here we need~$6\epsilon\rho\le1$). The various unions then contribute a linear factor in~$n$, which is absorbed into the exponential once~$n$ is sufficiently large.
\end{proofsect}

It is possible that a proper merge of the arguments in the previous two proofs might yield the same result without relying on Antal and Pisztora's bound \eqref{4.2}. (Indeed, the main other ``external'' ingredient of our proofs is Grimmett and Mar\-strand's paper~\cite{Grimmett-Marstrand} which lies at the core of~\cite{Antal-Pisztora} as well.) However, we find the argument using \eqref{4.2} conceptually cleaner and so we are content with the present, even though not necessarily optimal, proof.

\smallskip
Next we state a trivial, but interesting technical lemma:

\begin{mylemma}
\label{claim:caucschwa}
Let~$p>1$ and~$r\in[1,p)$. Suppose that~$X_1,X_2,\dots$ are random variables such that~$\sup_{j\ge1}\Vert X_j\Vert_p<\infty$ and let~$N$ be a random variable taking values in positive integers such that~$N\in L^s$ for some~$s$ satisfying
\begin{equation}
\label{E:s-def}
s>r\frac{1+\ffrac1p}{1-\ffrac rp}.
\end{equation}
Then $\sum_{j=1}^NX_j\in L^r$. Explicitly,
\begin{equation}
\Bigl\Vert\sum_{j=1}^NX_j\Bigr\Vert_r\le C\bigl(\,\sup_{j\ge1}\Vert X_j\Vert_p\bigr)
\bigl(\Vert N\Vert_s\bigr)^{s{\textstyle[}\ffrac1r-\ffrac1p{\textstyle]}},
\end{equation}
where $C$ is a finite constant depending only on~$p$, $r$ and~$s$.
\end{mylemma}

\begin{proofsect}{Proof}
Let us define $q\in(1,\infty)$ by~$\ffrac rp+\ffrac1q=1$. From the H\"older inequality and the uniform bound on~$\Vert X_j\Vert_p$ we get
\begin{equation}
\begin{aligned}
E\,\Bigl|\sum_{j=1}^NX_j\Bigr|^r
&=
\sum_{n\ge1}E\biggl(\,\Bigl|\sum_{j=1}^nX_j\Bigr|^r\1_{\{N=n\}}\biggr)
\\
&\le\sum_{n\ge1}\,\Bigl\Vert\sum_{j=1}^n X_j\Bigr\Vert_p^r\,P(N=n)^{\ffrac1q}
\\
&\le\bigl(\,\sup_{j\ge1}\Vert X_j\Vert_p\bigr)^r\,\sum_{n\ge1}n^r\,P(N=n)^{\ffrac1q}.
\end{aligned}
\end{equation}
Under the assumption that~$N$ has~$s$ moments, we get
\begin{equation}
\sum_{n\ge1}n^r\,P(N=n)^{\ffrac1q}
\le\biggl(\,\sum_{n\ge1}n^{(r-\ffrac sq)\ffrac pr}\biggr)^{\ffrac rp}\bigl(E(N^s)\bigr)^{\ffrac1q}
\end{equation}
by invoking the H\"older inequality one more time. The first term on the right-hand side is finite whenever~$s$ obeys the bound~\eqref{E:s-def}.
\end{proofsect}

\begin{proofsect}{Proof of Proposition~\ref{lemma:next_corrector}}
Let~$\chi(x,\omega)$ be the corrector. By Theorem~\ref{thm2.1}, on the set~$\{x\in\scrC_\infty\}$, $\chi(x,\cdot)$ is an $L^2$-limit of functions~$\chi_\epsilon(x,\cdot)=\psi_\epsilon\circ\tau_x-\psi_\epsilon$, as~$\epsilon\downarrow0$. To prove that~$\chi(v_e,\cdot)\in L^1$, recall the notation~$G_e^{(\epsilon)}$ from Lemma~\ref{lemma2.3} and let---as in Lemma~\ref{claim:pathlen}---$L=L(\omega)$ be the length of the shortest path from~$0$ to~$v_e$. Then
\begin{equation}
|\chi_\epsilon(v_e,\omega)|\le\sum_{x\colon|x|_\infty\le L(\omega)}\,\sum_{e\colon|e|=1}\bigl|G_e^{(\epsilon)}\circ\tau_x(\omega)\bigr|.
\end{equation}
But Theorem~\ref{thm2.1} ensures that~$\Vert G_e^{(\epsilon)}\circ\tau_x\Vert_2\le \Vert G_e^{(\epsilon)}\Vert_2<C$ for all~$x$ and~$e$ and all~$\epsilon>0$, while the number of terms in the sum does not exceed $N(\omega)=2d(2L(\omega)+1)^d$. By Lemma~\ref{claim:pathlen},~$N$ has all moments and so, by Lemma~\ref{claim:caucschwa}, $\sup_{\epsilon>0}\Vert\chi_\epsilon(v_e,\cdot)\Vert_r<\infty$ for all~$r\in[1,2)$. In particular,~$\chi(v_e,\cdot)\in L^1$.

In order to prove part~(2), we first note that a uniform bound on~$L^r$-norm of~$\chi_\epsilon(v_e,\cdot)$ for some~$r>1$ implies that the family~$\{\chi_\epsilon(v_e,\cdot)\}_{\epsilon>0}$ is uniformly integrable. Since $\chi_\epsilon(v_e,\cdot)\to\chi(v_e,\cdot)$ in probability, $\chi_\epsilon(v_e,\cdot)\to\chi(v_e,\cdot)$ in~$L^1$ and it thus suffices to prove
\begin{equation}
\E_0\bigl(\chi_\epsilon(v_e,\cdot)\bigr)=0,\qquad\epsilon>0.
\end{equation}
This is implied by Theorem~\ref{thm3.2} and the fact~$\chi_\epsilon(v_e,\cdot)=\psi_\epsilon\circ\sigma_e-\psi_\epsilon$ with~$\psi_\epsilon$ absolutely integrable.
\end{proofsect}

\begin{proofsect}{Proof of Theorem~\ref{thm4.1}}
Let $f(\omega)=\chi(n_1(\omega)e,\omega)$, and let $\sigma_e$ be the induced shift in the direction of~$e$.
Then we can write
\begin{equation}
\chi\bigl(n_k(\omega)e,\omega\bigr)=\sum_{\ell=0}^{k-1}f\circ\sigma_e^{\,\ell}(\omega).
\end{equation}
By Proposition~\ref{lemma:next_corrector} we have~$f\in L^1$ and $\E_0(f)=0$. Since Theorem~\ref{thm3.2} ensures that~$\sigma_e$ is~$\BbbP_0$-preserving and ergodic, the claim follows from Birkhoff's Ergodic Theorem.
\end{proofsect}

\section{Sublinearity everywhere}
\label{sec:sublin}\noindent
Here we will prove the principal technical estimates of this work. The level of control is different in~$d=2$ and~$d\ge3$, so we treat these cases separately. (Notwithstanding, the~$d\ge3$ proof applies in~$d=2$ as well.)

\subsection{Sublinearity in two dimensions}
We begin with an estimate of the corrector in large boxes in~$\Z^2$:

\begin{theorem}
\label{thm5.1}
Let~$d=2$ and let~$\chi$ be the corrector defined in Theorem~\ref{thm2.1}. Then for~$\BbbP_0$-almost every~$\omega\in\Omega_0$,
\begin{equation}
\label{5.1}
\lim_{n\to\infty}\,\max_{\begin{subarray}{c}
x\in\scrC_\infty(\omega)\\|x|_\infty\le n
\end{subarray}}
\frac{|\chi(x,\omega)|}n=0.
\end{equation}
\end{theorem}

The proof will be based on the following concept:

\begin{mydefinition}
\label{def5.6}
Given~$K>0$ and~$\epsilon>0$, we say that a site~$x\in\Z^d$ is $K,\epsilon$-\emph{good} (or just \emph{good}) in configuration~$\omega\in\Omega$ if $x\in\scrC_\infty(\omega)$ and 
\begin{equation}
\bigl|\chi(y,\omega)-\chi(x,\omega)\bigr|< K+\epsilon|x-y|
\end{equation}
holds for every~$y\in\scrC_\infty(\omega)$ of the 
form~$y=\ell e$, where $\ell\in\Z$ and~$e$ is a unit coordinate vector. 
We will use~$\scrG_{K,\epsilon}=\scrG_{K,\epsilon}(\omega)$ to denote the set of $K,\epsilon$-good sites in configuration~$\omega$.
\end{mydefinition}

On the basis of Theorem~\ref{thm4.1} it is clear that for each~$\epsilon>0$ there exists a~$K<\infty$ such that the $\BbbP_0(0\in\scrG_{K,\epsilon})>0$. Our first goal is to estimate the size of the largest interval free of good points in blocks~$[-n,n]$ on the coordinate axes:

\begin{mylemma}
\label{lemma5.3d}
Let~$e$ be one of the principal lattice vectors in~$\Z^2$ and, given~$\epsilon>0$, let~$K$ be so large that~$\BbbP_0(0\in\scrG_{K,\epsilon})>0$. For all~$n\ge1$ and~$\omega\in\Omega$, let~$y_0<\dots<y_r$ be the ordered set of all integers from~$[-n,n]$ such that~$y_ie\in\scrG_{K,\epsilon}(\omega)$. Let
\begin{equation}
\label{5.3}
\triangle_n(\omega)=\max_{j=1,\dots,r}(y_j-y_{j-1}).
\end{equation}
(If no such~$y_i$ exists, we define~$\triangle_n(\omega)=\infty$.) Then
\begin{equation}
\lim_{n\to\infty}\frac{\triangle_n}n=0,\qquad\BbbP\text{\rm-almost surely}.
\end{equation}
\end{mylemma}

\begin{proofsect}{Proof}
Since~$\BbbP$ is~$\tau_e$ invariant and~$\tau_e$ is ergodic, we have
\begin{equation}
\lim_{n\to\infty}\frac1{n+1}\sum_{k=0}^n\1_{\{0\in\scrG_{K,\epsilon}\}}\circ\tau_e^k=\BbbP(0\in\scrG_{K,\epsilon})
\end{equation}
$\BbbP$-almost surely. A similar statement applies to the limit~$n\to-\infty$. But if~$\triangle_n/n$ does not tend to zero, at least one of these limits would not exist.
\end{proofsect}

\begin{figure}[t]
\centerline{\epsfig{figure=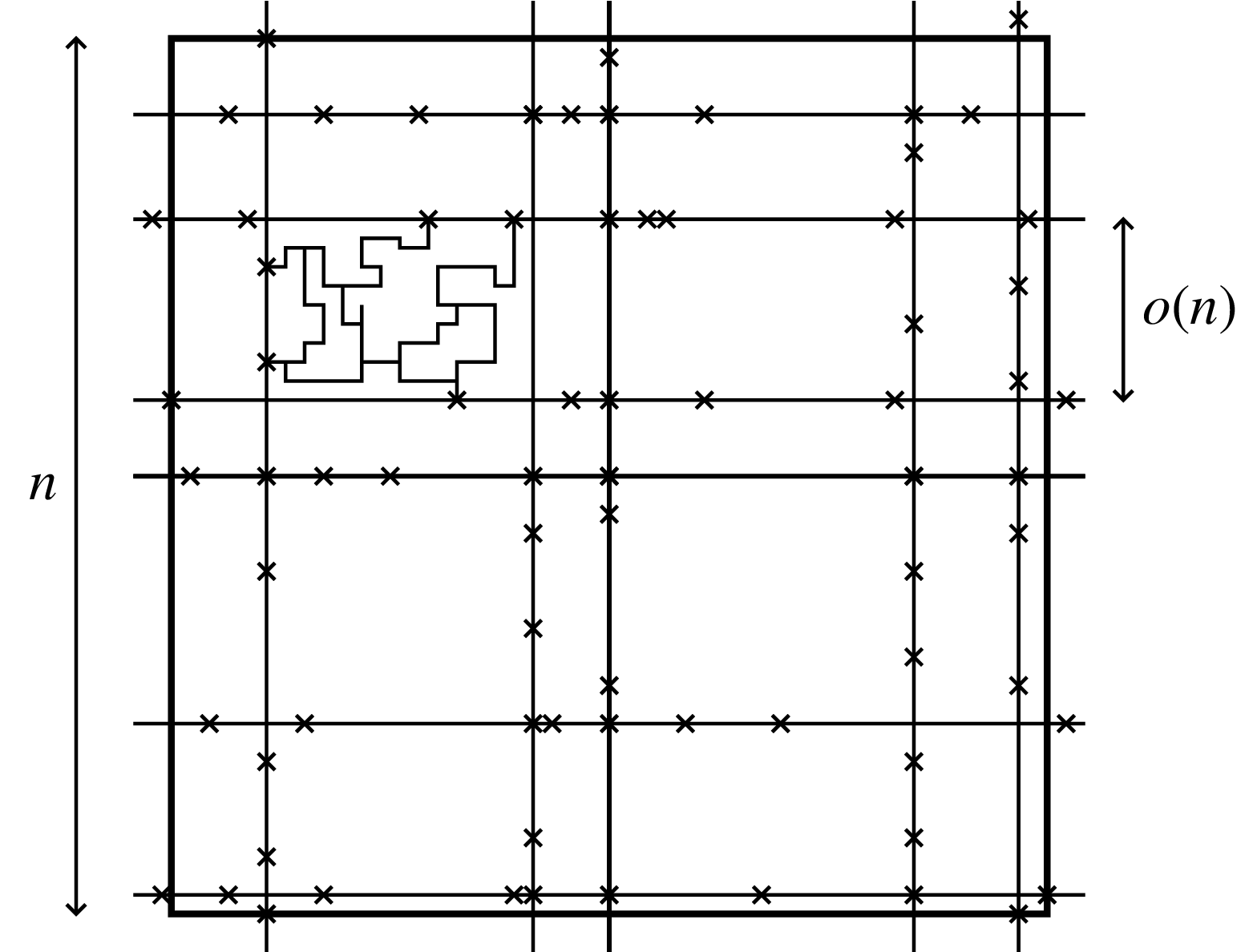, width=0.75\textwidth}}
\caption{An illustration of the main idea of the proof of Theorem~\ref{thm5.1}.
Here a square of side~$n$ is intersected by a grid~$\G$ of \emph{good lines} ``emanating'' from the good points on the~$x$ and~$y$ axes. The crosses represent the points on these lines which are in~$\scrC_\infty$. Along the good lines the corrector grows slower than linear and so anywhere on~$\G$ sublinearity holds. For the part of~$\scrC_\infty$ that is not on~$\G$, the maximum principle for~$x\mapsto x+\chi(x,\omega)$ lets us bound the corrector by the values on the parts of the grid that surround it, modulo factors of order~$o(n)$.}
\label{fig3}
\end{figure}  

\begin{proofsect}{Proof of Theorem~\ref{thm5.1}}
Fix~$\epsilon\in(0,\ffrac12)$ and let~$K_0$ be such that~$\BbbP(0\in\scrG_{K,\epsilon})>0$ for all~$K\ge K_0$ (we are using that~$\scrG_{K,\epsilon}$ increases with~$K$). Let~$\Omega_0^\star\subset\Omega_0$ be the set of configurations such that the conclusion of Lemma~\ref{lemma5.3d} applies for both~$x$ and~$y$-axes, 
and that shift-invariance \eqref{2.12d} holds for all~$x,y$ in the infinite cluster.
We will show that for every~$\omega\in\Omega_0^\star$ the limsup in~\eqref{5.1} is less than~$6\epsilon$ almost surely.

Let~$e_1$ and~$e_2$ denote the coordinate vectors in~$\Z^2$. Fix~$\omega\in\Omega_0^\star$ and adjust~$K\ge K_0$ so that $0\in\scrG_{K,\epsilon}$. (This is possible by the definition of~$\Omega_0^\star$.) 
Then we define~$(x_k)_{k\in\Z}$ to be the increasing two-sided sequence of all integers such that
$x_ke_1$ exhausts all $K,\epsilon$-good points on the~$e_1$-axis, i.e.,
\begin{equation}
\label{5.6}
x_ke_1\in\scrG_{K,\epsilon}(\omega),\qquad k\in\Z.
\end{equation}
If~$\triangle_n$ be the maximal gap between consecutive~$x_j$'s that lie in~$[-n,n]$, cf~\eqref{5.3}, we define~$n_1(\omega)$ be the least integer such that~$\triangle_n/n<\epsilon$ for all~$n\ge n_1(\omega)$. Similarly we identify a two-sided 
increasing sequence~$(y_n)_{n\in\Z}$ of integers exhausting the sites such that
\begin{equation}
\label{5.7}
y_ke_2\in\scrG_{K,\epsilon}(\omega),\qquad k\in\Z,
\end{equation}
and let~$n_2(\omega)$ be the quantity corresponding to~$n_1(\omega)$ in this case.

Let~$n_0=\max\{n_1,n_2\}$. We claim that for all~$n\ge n_0(\omega)$,
\begin{equation}
\label{5.8}
\max_{\begin{subarray}{c}
x\in\scrC_\infty(\omega)\\|x|_\infty\le n
\end{subarray}}
|\chi(x,\omega)|\le 2K+6\epsilon n.
\end{equation}
To prove this, let us consider the grid~$\G=\G(\omega)$ of \emph{good lines}
\begin{equation}
\{x_ke_1+ne_2\colon n\in\Z\},\qquad k\in\Z,
\end{equation}
and
\begin{equation}
\{ne_1+y_ke_2\colon n\in\Z\},\qquad k\in\Z,
\end{equation}
see Fig.~\ref{fig3}.
As a first step we will use the harmonicity of~$x\mapsto x+\chi(x,\omega)$ to deal with~$x\in\scrC_\infty\setminus\G$. Indeed, any such~$x$ is enclosed between two horizontal and two vertical grid lines and every path on~$\scrC_\infty$ connecting~$x$ to ``infinity'' necessarily intersects one of these lines at a point which is also in~$\scrC_\infty$. Applying the maximum (and minimum) principle for harmonic functions we get
\begin{equation}
\label{5.11}
\max_{\begin{subarray}{c}
x\in\scrC_\infty\smallsetminus\G\\|x|_\infty\le n
\end{subarray}}
|\chi(x,\omega)|
\le 2\epsilon n+
\max_{\begin{subarray}{c}
x\in\scrC_\infty\cap\G\\|x|_\infty\le 2n
\end{subarray}}
|\chi(x,\omega)|.
\end{equation}
Here we used that the enclosing lines are not more than $\frac\epsilon{1-\epsilon}n\le 2\epsilon n\le n$ apart and, in particular, they all intersect the block $[-2n,2n]\times[-2n,2n]$.

To estimate the maximum on the grid, we pick, say, a horizontal grid line with~$y$-coordinate~$y_k$ and note that, by~\eqref{2.12d}, for every~$x\in\scrC_\infty$ on this line,
\begin{equation}
\chi(x,\omega)-\chi(y_ke_2,\omega)=\chi(x-y_ke_2,\tau_{y_ke_2}\omega).
\end{equation}
By~\eqref{5.7} and the fact that~$x-y_ke_2\in\scrC_\infty(\tau_{y_ke_2}\omega)$ we have
\begin{equation}
\bigl|\chi(x,\omega)-\chi(y_ke_2,\omega)\bigr|\le K+2\epsilon n
\end{equation}
whenever~$x$ is such that~$|x|_\infty\le2n$. Applying the same argument to the vertical line through the origin, and~$x$ replaced by~$y_ke_2$, we get
\begin{equation}
\bigl|\chi(x,\omega)\bigr|\le 2K+4\epsilon n
\end{equation}
for every~$x\in\scrC_\infty\cap\G$ with~$|x|_\infty\le 2n$. Combining this with~\eqref{5.11}, the estimate~\eqref{5.8} and the whole claim are finally proved.
\end{proofsect}

Interestingly, a variant of the above strategy for controlling the corrector in $d=2$ has independently been developed by Chris Hoffman~\cite{Hoffman-geodesics} to control the geodesics in the first-passage percolation on~$\Z^2$.

\subsection{Three and higher dimensions}
In~$d\ge3$ we have the following weaker version of Theorem~\ref{thm5.1}:

\begin{theorem}
\label{thm5.4}
Let~$d\ge3$. Then for all~$\epsilon>0$ and~$\BbbP_0$-almost all~$\omega$,
\begin{equation}
\label{5.15c}
\limsup_{n\to\infty}\frac1{(2n+1)^d}\sum_{\begin{subarray}{c}
x\in\scrC_\infty(\omega)\\|x|\le n
\end{subarray}}
\1_{\{|\chi(x,\omega)|\ge\epsilon n\}}=0.
\end{equation}
\end{theorem}

Here we fix the dimension~$d$ and run an induction over $\nu$-dimensional sections of the $d$-dimensional box~$\{x\in\Z^d\colon|x|\le n\}$. Specifically, for each~$\nu=1,\dots,d$, let~$\Lambda_n^\nu$ be the~$\nu$-dimensional box
\begin{equation}
\Lambda_n^\nu=\bigl\{k_1e_1+\dots+k_\nu e_\nu\colon k_i\in\Z,\,|k_i|\le n\,\,\forall i=1,\dots,\nu\bigr\}.
\end{equation}
The induction eventually gives~\eqref{5.15c} for~$\nu=d$ thus proving the theorem.

Since it is not advantageous to assume that~$0\in\scrC_\infty$, we will carry out the proof for \emph{differences} of the form $\chi(x,\omega)-\chi(y,\omega)$ with~$x,y\in\scrC_\infty$.
For each~$\omega\in\Omega$, we thus consider the (upper) density
\begin{equation}
\label{5.17}
\varrho_\nu(\omega)=
\lim_{\epsilon\downarrow0}\,
\limsup_{n\to\infty}
\inf_{y\in\scrC_\infty(\omega)\cap\Lambda_n^1}\,
\frac1{|\Lambda_n^\nu|}\,
\sum_{x\in\scrC_\infty(\omega)\cap\Lambda_n^\nu}
\1_{\{|\chi(x,\omega)-\chi(y,\omega)|\ge\epsilon n\}}.
\end{equation}
Note that the infimum is taken only over sites in one-dimensional box $\Lambda_n^1$.
Our goal is to show by induction that~$\varrho_\nu=0$ almost surely for all~$\nu=1,\dots,d$. The induction step is encapsulated into the following lemma:

\begin{mylemma}
\label{lemma5.5}
Let~$1\le\nu<d$. If $\varrho_\nu=0$,~$\BbbP$-almost surely, then also~$\varrho_{\nu+1}=0$, $\BbbP$-almost surely.
\end{mylemma}

Before we start the formal proof, let us discuss its main idea: 
Suppose that $\varrho_\nu=0$ for some~$\nu<d$, $\BbbP$-almost surely.
Pick $\epsilon>0$. Then for $\BbbP$-almost every~$\omega$ and all sufficiently large~$n$, there exists a set of sites~$\Delta\subset\Lambda_n^\nu\cap\scrC_\infty$ such that
\begin{equation}
\label{5.19}
\bigl|(\Lambda_n^\nu\cap\scrC_\infty)\setminus\Delta\bigr|
\le \epsilon|\Lambda_n^\nu|
\end{equation}
and
\begin{equation}
\label{5.19k}
\bigl|\chi(x,\omega)-\chi(y,\omega)\bigr|\le\epsilon n,\qquad x,y\in\Delta.
\end{equation}
Moreover,~$n$ sufficiently large,~$\Delta$ could be picked so that $\Delta\cap\Lambda_n^1\ne\emptyset$ and, assuming~$K\gg1$, the non-$K,\epsilon$-good sites could be pitched out with little loss of density to achieve even
\begin{equation}
\label{5.21}
\Delta\subset\scrG_{K,\epsilon}.
\end{equation}
(All these claims are direct consequences of the Pointwise Ergodic Theorem and the fact that~$\BbbP(0\in\scrG_{K,\epsilon})$ converges to the density of~$\scrC_\infty$ as~$K\to\infty$.)

\begin{figure}[t]
\centerline{\epsfig{figure=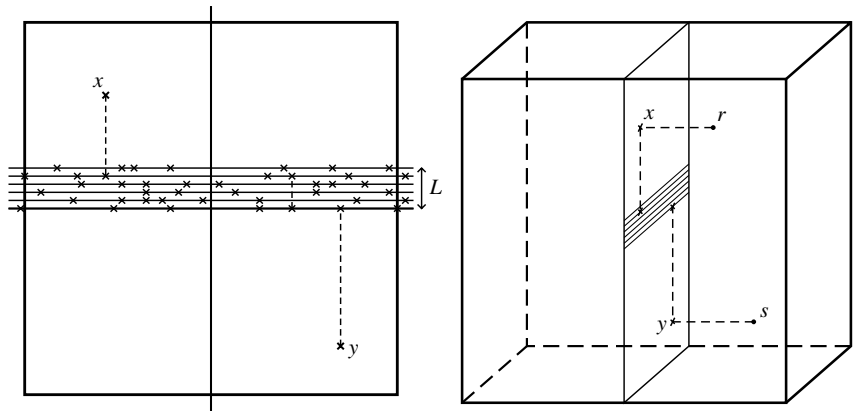, width=\textwidth}}
\caption{The main idea underlying the proof of Theorem~\ref{thm5.4}.
The figure on the left represents an $n\times n$ square in a two-dimensional plane in~$\Z^3$; the crosses now stand for \emph{good} sites; cf Definition~\ref{def5.6}. Here~$L$ is chosen so that~$(1-\delta)$-fraction of all vertical lines find a good point on the intersection with one of the~$L$ horizontal lines;~$n$ is assumed so large that every pair of these lines has two good points ``above'' each other. Any two good points~$x$ and~$y$ in the square are connected by broken-line path that uses at most 4 good points in between. The dashed lines indicate the vertical pieces of one such path. The figure on the right indicates how this is used to control the difference of the corrector at two general points~$r,s\in\scrC_\infty$ in an~$n\times n\times n$ cube in~$\Z^3$---with obvious extensions to all~$d\ge3$.}
\label{fig4}
\end{figure}

As a result of this construction we have
\begin{equation}
\label{5.22}
\bigl|\chi(z,\omega)-\chi(x,\omega)\bigr|\le K+\epsilon n
\end{equation}
for any~$x\in\Delta$ and any~$z\in\Lambda_n^{\nu+1}\cap\scrC_\infty$ of the form~$x+je_{\nu+1}$. Thus, if~$r,s\in\scrC_\infty\cap\Lambda_n^{\nu+1}$ are of the latter form, $r=x+je_{\nu+1}$ and~$s=y+ke_{\nu+1}$---see Fig.~\ref{fig4} for an illustration---then \eqref{5.22} implies
\begin{equation}
\label{5.22ba}
|\chi(r,\omega)-\chi(s,\omega)|\le |\chi(x,\omega)-\chi(y,\omega)|+2K+2\epsilon n.
\end{equation}
Invoking the ``induction hypothesis'' \eqref{5.19k}, the right-hand side is less than $2K+3\epsilon n$, implying a bound of the type \eqref{5.19k} but one-dimension higher.

Unfortunately, the above is not sufficient to prove \eqref{5.19k} for all but a vanishing fraction of all sites in~$\Lambda_n^{\nu+1}$. The reason is that the~$r$'s and~$s$'s for which \eqref{5.22ba} holds need to be of the form~$x+je_{\nu+1}$ for some~$x\in\Delta\cap\scrC_\infty$. But~$\scrC_\infty$ will occupy only about $P_\infty=\BbbP(0\in\scrC_\infty)$ fraction of all sites in~$\Lambda_n^\nu$, and so this argument does not permit us to control more than fraction about~$P_\infty$ of~$\Lambda_n^{\nu+1}\cap\scrC_\infty$.

To fix this problem, we will have to work with a ``stack'' of translates of~$\Lambda_n^\nu$ at the same time. (These correspond to the stack of horizontal lines on the left of of Fig.~\ref{fig4}.) Explicitly, consider the collection of~$\nu$-boxes
\begin{equation}
\Lambda_{n,j}^\nu=\tau_{je_{\nu+1}}(\Lambda_n^\nu),\qquad j=1,\dots,L.
\end{equation}
Here~$L$ is a deterministic number chosen so that, for a given~$\delta>0$, the set
\begin{equation}
\Delta_0=\bigl\{x\in\Lambda_n^\nu\colon\exists j\in\{0,\dots,L-1\},\, x+je_{\nu+1}\in\Lambda_{n,j}^\nu\cap\scrC_\infty\bigr\}
\end{equation}
is so large that
\begin{equation}
|\Delta_0|\ge(1-\delta)|\Lambda_n^\nu|
\end{equation}
once~$n$ is sufficiently large. These choices ensure that~$(1-\delta)$-fraction of~$\Lambda_n^\nu$ is now ``covered'' which by repeating the above argument gives us control over~$\chi(r,\omega)$ for nearly the same fraction of all sites~$r\in\Lambda_n^{\nu+1}\cap\scrC_\infty$.

\begin{proofsect}{Proof of Lemma~\ref{lemma5.5}}
Let~$\nu<d$ and suppose that $\varrho_\nu=0$, $\BbbP$-almost surely.  
Fix~$\delta$ with~$0<\delta<\frac12 P_\infty^2$ and let~$L$ be as defined above. Choose~$\epsilon>0$ so that
\begin{equation}
\label{5.26}
L\epsilon+\delta<\frac12P_\infty^2.
\end{equation}
For a fixed but large~$K$, and $\BbbP$-almost every~$\omega$ and $n$ exceeding an $\omega$-dependent quantity, for each $j=1,\dots,L$, we can find~$\Delta_j\subset\Lambda_{n,j}^\nu\cap\scrC_\infty$ satisfying the properties \twoeqref{5.19}{5.21}---with~$\Lambda_n^\nu$ replaced by $\Lambda_{n,j}^\nu$.
Given $\Delta_1,\dots,\Delta_L$, let~$\Lambda$ be the set of sites in~$\Lambda_n^{\nu+1}\cap\scrC_\infty$ whose projection onto the linear subspace $\BbbH=\{k_1e_1+\dots+k_\nu e_\nu\colon k_i\in\Z\}$ belongs to the corresponding projection of $\Delta_1\cup\dots\cup\Delta_L$.
Note that the~$\Delta_j$ could be chosen so that~$\Lambda\cap\Lambda_n^1\ne\emptyset$.

By their construction, the projections of the~$\Delta_j$'s, $j=1,\dots,L$, onto $\BbbH$ ``fail to cover'' at most~$L\epsilon|\Lambda_n^\nu|$ sites in~$\Delta_0$, and so at most~$(\delta+L\epsilon)|\Lambda_n^\nu|$ sites in~$\Lambda_n^\nu$ are not of the form~$x+i e_{\nu+1}$ for some $x\in\bigcup_j\Delta_j$. It follows that
\begin{equation}
\label{5.21a}
\bigl|(\Lambda_n^{\nu+1}\cap\scrC_\infty)\setminus\Lambda\bigr|
\le(\delta+L\epsilon)|\Lambda_n^{\nu+1}|,
\end{equation}
i.e.,~$\Lambda$ contains all except at most~$(L\epsilon+\delta)$-fraction of all sites in~$\Lambda_n^{\nu+1}$ that we care about.
Next we note that if~$K$ is sufficiently large, then for every~$1\le i<j\le L$, the set~$\BbbH$ contains at least~$\frac12 P_\infty^2$-fraction of sites~$x$ such that
\begin{equation}
\label{5.27}
z_i\,\overset{\text{\rm def}}=\,x+ie_\nu\in\scrG_{K,\epsilon}
\quad\text{and}\quad
z_j\,\overset{\text{\rm def}}=\,x+je_\nu\in\scrG_{K,\epsilon}.
\end{equation}
Since we assumed~\eqref{5.26}, once~$n\gg1$, for each pair~$(i,j)$ with $1\le i<j\le L$ such~$z_i$ and~$z_j$ can be found so that~$z_i\in\Delta_i$ and~$z_j\in\Delta_j$. But the~$\Delta_j$'s were picked to make \eqref{5.19k} true and so via these pairs of sites we now show that
\begin{equation}
\label{5.28}
\bigl|\chi(y,\omega)-\chi(x,\omega)\bigr|\le K+\epsilon L+2\epsilon n
\end{equation}
for every~$x,y\in\Delta_1\cup\dots\cup\Delta_L$; see again (the left part of) Fig.~\ref{fig4}.

From \eqref{5.19k} and \eqref{5.28} we now conclude that for all~$r,s\in\Lambda$,
\begin{equation}
\label{5.22a}
\bigl|\chi(r,\omega)-\chi(s,\omega)\bigr|\le 3K+\epsilon L+4\epsilon n<5\epsilon n,
\end{equation}
provided that~$\epsilon n>3K+\epsilon L$. If~$\varrho_{\nu,\epsilon}$ denotes the right-hand side of \eqref{5.17} before taking~$\epsilon\downarrow0$, the bounds \eqref{5.21a} and \eqref{5.22a} and $\Lambda\cap\Lambda_n^1\ne\emptyset$ yield
\begin{equation}
\varrho_{\nu+1,5\epsilon}(\omega)\le\delta+L\epsilon,
\end{equation}
for~$\BbbP$-almost every~$\omega$. But the left-hand side of this inequality increases as~$\epsilon\downarrow0$ while the right-hand side decreases. Thus, taking~$\epsilon\downarrow0$ and~$\delta\downarrow0$ proves that~$\rho_{\nu+1}=0$ holds~$\BbbP$-almost surely.
\end{proofsect}

\begin{proofsect}{Proof of Theorem~\ref{thm5.4}}
The proof is an easy consequence of Lemma~\ref{lemma5.5}.
First, by Theorem~\ref{thm4.1} we know that $\varrho_1(\omega)=0$ for $\BbbP_0$-almost every~$\omega$. 
Invoking appropriate shifts, the same conclusion applies~$\BbbP$-almost surely. Using induction on dimension, Lemma~\ref{lemma5.5} then tells us that $\varrho_d(\omega)=0$ for~$\BbbP_0$-almost every~$\omega$. Let~$\omega\in\Omega_0$. By Theorem~\ref{thm4.1}, for each~$\epsilon>0$ there is $n_0=n_0(\omega)$ with $\BbbP_0(n_0<\infty)=1$ such that for all~$n\ge n_0(\omega)$, we have $|\chi(x,\omega)|\le\epsilon n$ for all $x\in\Lambda_n^1\cap\scrC_\infty(\omega)$.
Using this to estimate away the infimum in~\eqref{5.17}, the fact that~$\varrho_d=0$ now immediately implies~\eqref{5.15c} for all~$\epsilon>0$.
\end{proofsect}

\section{Proof of main results}
\label{sec:invprnp}\noindent
Here we will finally prove our main theorems. First, in Sect.~\ref{sec6.1}, we will show the convergence of the ``lazy'' walk on the deformed graph to Brownian motion and then, in Sect.~\ref{sec6.2}, we use our previous results on corrector growth to extend this to the walk on the original graph. This separation will allow us to treat the parts of the proof common for~$d=2$ and~$d\ge3$ in a unified way. Theorem~\ref{thm:2ndmainthm}, which concerns the ``agile'' walk, is proved in Sect.~\ref{sec6.3}.

\subsection{Convergence on deformed graph}
\label{sec6.1}
We begin with a simple observation that will drive all underlying derivations:

\begin{mylemma}
\label{lemma6.2}
Fix~$\omega\in\Omega_0$ and let~$x\mapsto\chi(x,\omega)$ be the corrector. Given a path of random walk~$(X_n)_{n\ge0}$ with law~$P_{0,\omega}$, let
\begin{equation}
\label{E:6.1}
M_n^{(\omega)}=X_n+\chi(X_n,\omega),\qquad n\ge0.
\end{equation}
Then~$(M_n^{(\omega)})_{n\ge0}$ is an~$L^2$-martingale for the filtration~$(\sigma(X_0,\dots,X_n))_{n\ge0}$.
Moreover, conditional on~$X_{k_0}=x$, the increments $(M_{k+k_0}^{(\omega)}-M_{k_0}^{(\omega)})_{k\ge0}$ have the same law as~$(M_k^{(\tau_x\omega)})_{k\ge0}$.
\end{mylemma}

\begin{proofsect}{Proof}
Since~$X_n$ is bounded, $\chi(X_n,\omega)$ is bounded and so~$M_n^{(\omega)}$ is square integrable with respect to~$P_{0,\omega}$. Since~$x\mapsto x+\chi(x,\omega)$ is harmonic with respect to the transition probabilities of the random walk~$(X_n)$ with law~$P_{0,\omega}$, we have
\begin{equation}
E_{0,\omega}\bigl(M_{n+1}^{(\omega)}\big|\sigma(X_n)\bigr)=M_n^{(\omega)},\qquad n\ge0,
\end{equation}
$P_{0,\omega}$-almost surely. Since~$M_n^{(\omega)}$ is~$\sigma(X_n)$-measurable,~$(M_n^{(\omega)})$ is a martingale. 
The stated relation between the laws of~$(M_{k+k_0}^{(\omega)}-M_{k_0}^{(\omega)})_{k\ge0}$ and~$(M_k^{(\tau_x\omega)})_{k\ge0}$ is implied by the shift-invariance \eqref{2.12d} and the fact that~$(M_n^{(\omega)})$ is a simple random walk on the deformed infinite component.
\end{proofsect}

Next we will establish the convergence of the above martingale to Brownian motion.
The precise statement is as follows:

\begin{theorem}
\label{thm6.1}
Let~$d\ge2$, $p>p_\cc$ and~$\omega\in\Omega_0$. Let~$(X_n)_{n\ge0}$ be the random walk with law~$P_{0,\omega}$ and let $(M_n^{(\omega)})_{n\ge0}$ be as defined in~\eqref{E:6.1}. Let $(\widehat B_n^{(\omega)}(t)\colon t\ge0)$ be defined by
\begin{equation}
\label{E:1.5a}
\widehat B_n^{(\omega)}(t)=\frac1{\sqrt n}\bigl(
M_{\lfloor tn\rfloor}^{(\omega)}+(tn-\lfloor tn\rfloor)(M_{\lfloor tn\rfloor+1}^{(\omega)}-M_{\lfloor tn\rfloor}^{(\omega)})\bigr),
\qquad t\ge0.
\end{equation}
Then for all~$T>0$ and~$\BbbP_0$-almost every~$\omega$, the law of~$(\widehat B_n(t)\colon 0\le t\le T)$ on $(C[0,T],\scrW_T)$ converges weakly to the law of an isotropic Brownian motion $(B_t\colon 0\le t\le T)$ with diffusion constant~$D$, i.e.,~$E(B_t^2)=Dt$, where
\begin{equation}
\label{E:D-eq}
D=\E_0\Bigl(E_{0,\omega}\bigl|X_1+\chi(X_1,\omega)\bigr|^2\Bigr)\in(0,\infty).
\end{equation}
\end{theorem}

\begin{proofsect}{Proof}
Without much loss of generality, we may confine ourselves to the case when~$T=1$.
Let~$\scrF_k=\sigma(X_0,X_1,\dots,X_k)$ and fix a vector~$a\in\R^d$. We will show that (the piece-wise linearization) of~$t\mapsto a\cdot M_{\lfloor tn\rfloor}^{(\omega)}$ scales to one-dimensional Brownian motion. For~$m\le n$, consider the random variable
\begin{equation}
\label{6.4a}
V_{n,m}^{(\omega)}(\epsilon)=\frac1n\sum_{k=0}^{m}E_{0,\omega}\Bigl(\bigl[a\cdot(M_{k+1}^{(\omega)}-M_k^{(\omega)})\bigr]^2\1_{\{|a\cdot(M_{k+1}^{(\omega)}-M_k^{(\omega)})|\ge\epsilon\sqrt n\}}\Big|\scrF_k\Bigr).
\end{equation}
In order to apply the Lindeberg-Feller Functional~CLT for martingales
(Theorem~7.7.3 of Durrett~\cite{Durrett}), we need to verify that for~$\BbbP_0$-almost every~$\omega$,
\begin{enumerate}
\item[(1)]
$V_{n,\lfloor tn\rfloor}^{(\omega)}(0)\to Ct$ in~$P_{0,\omega}$-probability for all~$t\in[0,1]$ and some~$C\in(0,\infty)$.
\item[(2)]
$V_{n,n}^{(\omega)}(\epsilon)\to0$ in~$P_{0,\omega}$-probability for all~$\epsilon>0$.
\end{enumerate}
Both of these conditions will be implied by Theorem~\ref{thm3.1}. Indeed, by the last conclusion of Lemma~\ref{lemma6.2} we may write
\begin{equation}
V_{n,m}^{(\omega)}(\epsilon)=\frac1n\sum_{k=0}^{m}f_{\epsilon\sqrt n}\circ\tau_{X_k}(\omega),
\end{equation}
where
\begin{equation}
f_K(\omega)=E_{0,\omega}\Bigl(\bigl[a\cdot M_1^{(\omega)}\bigr]^2\1_{\{|a\cdot M_1^{(\omega)}|\ge K\}}\Bigr).
\end{equation}
Now if~$\epsilon=0$, Theorem~\ref{thm3.1} tells us that, for~$\BbbP_0$-almost every~$\omega$,
\begin{equation}
\lim_{n\to\infty}V_{n,n}^{(\omega)}(0)=\E_0\Bigl(E_{0,\omega}\bigl([a\cdot M_1^{(\omega)}]^2\bigr)\Bigr)=
\frac1d D|a|^2,
\end{equation}
where we used the symmetry of the joint expectations under rotations by~$90^\circ$. From here condition~(1) follows by scaling out the~$t$-dependence first and working with~$tn$ instead of~$n$.

On the other hand, when~$\epsilon>0$, we have~$f_{\epsilon\sqrt n}\le f_K$ once~$n$ is sufficiently large and so,~$\BbbP_0$-almost surely,
\begin{equation}
\label{6.8a}
\limsup_{n\to\infty}V_{n,n}^{(\omega)}(\epsilon)\le
\E_0\Bigl(E_{0,\omega}\bigl([a\cdot M_1^{(\omega)}]^2\1_{\{|a\cdot M_1^{(\omega)}|\ge K\}}\bigr)\Bigr)
\,\underset{K\to\infty}\longrightarrow\,0,
\end{equation}
where to apply Dominated Convergence we used that~$a\cdot M_1^{(\omega)}\in L^2$.
Hence, the above conditions~(1) and (2) hold---in fact, even with limits taken $P_{0,\omega}$-almost surely.
Applying the Martingale functional CLT and the Cram\'er-Wold device (Theorem~2.9.2 of \cite{Durrett}), we conclude that, for~$\BbbP_0$-almost every~$\omega$, the linear interpolation of the sequence $(M_k^{(\omega)}/\sqrt n)_{k=1,\dots,n}$ converges to isotropic Brownian motion with covariance matrix~$\frac1d D\1$.

To make the proof complete, we need to show that $D\in(0,\infty)$. Here the finiteness is immediate by the square-integrability of~$\chi$. The positivity can be shown in many ways: either by a direct computation from \eqref{E:D-eq} using that~$\E_0(E_{0,\omega}(X_1\cdot\chi(X_1,\omega))=0$ [which in turn is implied by $\E_0(\chi(e,\omega)\1_{\{\omega_e=1\}})=0$ for every coordinate vector~$e$] or by invoking the sublinearity of the corrector proved in Theorems~\ref{thm5.1}--\ref{thm5.4}, or by an appeal to the lower (or, alternatively, upper) bound in~\cite[Theorem~1]{Barlow}.
\end{proofsect}

\subsection{Correction on the corrector}
\label{sec6.2}\noindent
It remains to estimate the influence of the harmonic deformation on the path of the walk. As already mentioned, while our proof in~$d=2$ is completely self-contained, for~$d\ge3$ we rely heavily on (a discrete version of) the sophisticated Theorem~1 of Barlow~\cite{Barlow}.

\smallskip
Let us first dismiss the two-dimensional case of Theorem~\ref{thm:mainthm}:

\begin{proofsect}{Proof of Theorem~\ref{thm:mainthm} ($d=2$)}
We need to extend the conclusion of Theorem~\ref{thm6.1} to the linear interpolation of~$(X_n)$.
Since the corrector is an additive perturbation of~$M_n^{(\omega)}$, it clearly suffices to show that, for~$\BbbP_0$-almost every~$\omega$,
\begin{equation}
\max_{1\le k\le n}\frac{|\chi(X_k,\omega)|}{\sqrt n}\,\underset{n\to\infty}\longrightarrow\,0,\qquad \text{in $P_{0,\omega}$-probability}.
\end{equation}
By Theorem~\ref{thm5.1} we know that for every~$\epsilon>0$ there exists a~$K=K(\omega)<\infty$ such that
\begin{equation}
|\chi(x,\omega)|\le K+\epsilon|x|_\infty,\qquad x\in\scrC_\infty(\omega).
\end{equation}
If~$\epsilon<\ffrac12$, then this implies
\begin{equation}
|\chi(X_k,\omega)|\le2K+2\epsilon|M_k^{(\omega)}|.
\end{equation}
But the above CLT for~$(M_k)$ tells us that~$\max_{k\le n}|M_k^{(\omega)}|/\sqrt n$ converges in law to the maximum of a Brownian motion~$B(t)$ over~$t\in[0,1]$. Hence, if~$P$ denotes the probability law of the Brownian motion, the Portmanteau Theorem (Theorem~2.1 of~\cite{Billingsley}) allows us to conclude
\begin{equation}
\limsup_{n\to\infty}
P_{0,\omega}\bigl(\,\max_{k\le n}|\chi(X_k,\omega)|\ge\delta\sqrt n\,\bigr)
\le P\Bigl(\,\max_{0\le t\le1}|B(t)|\ge \frac{\delta}{2\epsilon}\Bigr).
\end{equation}
The right-hand side tends to zero as~$\epsilon\downarrow0$ for all~$\delta>0$.
\end{proofsect}

In order to prove the same result in~$d\ge3$, we will need the following upper bounds on the transition probability of our random walk:

\begin{theorem}
\label{thm6.3a}
(1)
There is a random variable~$C=C(\omega)$ with $\BbbP_0(C<\infty)=1$ such that for all~$\omega\in\Omega_0$ and all~$x\in\scrC_\infty(\omega)$,
\begin{equation}
\label{6.14q}
P_{0,\omega}(X_n=x)\le \frac{C(\omega)}{n^{d/2}},
\qquad n\ge1.
\end{equation}

\smallskip\noindent
(2)
There are constants~$\c_1,\c_2\in(0,\infty)$ and random
variables~$N_x=N_x(\omega)$ such that for all $\omega\in\Omega_0$, all~$x\in\scrC_\infty(\omega)$, 
all $R\ge1$, and all~$n\ge N_x(\omega)$,
\begin{equation}
P_{x,\omega}\bigl(|X_n-x|>R\bigr)\le 
\c_1\exp\bigl\{-\c_2R^2/n\bigr\}.
\end{equation}
Moreover, the random variables~$(N_x)$ have stretched-exponential tails, i.e., there exist constants~$\c_3>0$ and $\theta>0$ such that for all~$x\in\Z^d$,
\begin{equation}
\label{Nxtail}
\BbbP_0(N_x>R)\le \texte^{-\c_3R^\theta},\qquad R\ge1.
\end{equation}
\end{theorem}

For a continuous-time version of our walk, these bounds are the content of Theorem~1 of Barlow~\cite{Barlow}. 
(In fact, the continuous-time version of the bound \eqref{6.14q} was obtained already by Mathieu and Remy~\cite{Mathieu-Remy}.)
Unfortunately, to derive Theorem~\ref{thm6.3a} from Barlow's Theorem~1, one needs to invoke various non-trivial facts about percolation and/or mixing of Markov chains.
In Appendix~\ref{appA} we list these facts and show how to assemble all ingredients together to establish the above upper bounds.

\smallskip
\begin{proofsect}{Proof of Theorem~\ref{thm:mainthm} ($d\ge3$)}
We will adapt (the easier part of) the proof of Theorem~1.1 in Sidoravicius and Sznitman~\cite{Sidoravicius-Sznitman}.
First we show that the laws of $(\widetilde B_n(t)\colon t\le T)$ on $(C[0,T],\scrW_T)$ are tight. To that end it suffices to show (e.g., by Theorem~8.6 of Ethier-Kurtz~\cite{Ethier-Kurtz}) that if~$\scrS_n$ is the class of all stopping times of the filtration~$(\sigma(\{\widetilde B_n(s)\colon s\le t\}))_{0\le t\le T}$, then
\begin{equation}
\label{6.15a}
\limsup_{\epsilon\downarrow0}\,\limsup_{n\to\infty}\sup_{\tau\in\scrS_n}E_{0,\omega}
\bigl(|\widetilde B_n(\tau+\epsilon)-\widetilde B_n(\tau)|^2\bigr)=0.
\end{equation}
As in~\cite{Sidoravicius-Sznitman}, we replace~$\tau$ by its integer-valued approximation. Explicitly, let~$\hat\tau=\lfloor n\tau\rfloor+1$ and let~$\delta$ be a number such that~$n\delta=\lfloor n\epsilon\rfloor+1$. Since~$\hat\tau$ differs from~$n\tau$ by a constant of order unity, and similarly for~$\hat\tau+n\delta$ and~$n(\tau+\epsilon)$, we have
\begin{equation}
\label{E:6.18a}
|\widetilde B_n(\tau+\epsilon)-\widetilde B_n(\tau)|\le
\frac1{\sqrt n}|X_{\hat\tau+n\delta}-X_{\hat\tau}|+\frac{\c_4}{\sqrt n}
\end{equation}
for some constant~$\c_4<\infty$. This allows us to estimate \eqref{6.15a} by means of the second moment of $|X_{\hat\tau+n\delta}-X_{\hat\tau}|$.

Recalling that~$\tau\le T$, we may assume that~$\hat\tau\le 2Tn$. By~\eqref{Nxtail} we know that there exists an almost-surely finite random variable~$C'=C'(\omega)$ such that $\max_{|x|\le R}N_x\le C'(\omega)(\log R)^\zeta$ once~$R\ge2$, where~$\zeta=\ffrac2\theta$. Since~$|X_{\hat\tau}|\le 2Tn$, this implies that~$N_{X_{\hat\tau}}\le C'(\omega)[\log(2Tn)]^\zeta$. Theorem~\ref{thm6.3a}(2) and the strong Markov property---$\hat\tau$ is a stopping time of the random walk---tell us that, for some constant~$\c_5<\infty$ (depending only on~$\c_1,\c_2$ and the dimension),
\begin{equation}
\label{E:6.19a}
E_{0,\omega}\bigl(|X_{\hat\tau+n\delta}-X_{\hat\tau}|^2\bigr)\le \c_5\,\epsilon n,
\qquad n\ge n_0(\omega).
\end{equation}
Here we used~$\epsilon-\delta=O(\ffrac1n)$ and let~$n_0(\omega)$ be such that $\delta n\ge C'(\omega)[\log(2Tn)]^\zeta$ for all~$n\ge n_0$. The bound \eqref{6.15a} is now proved by combining \twoeqref{E:6.18a}{E:6.19a} and taking the required limits.

Once we know that the laws of~$(\widetilde B_n(t)\colon t\le T)$ are tight, it suffices to show the convergence of finite-dimensional distributions. In light of Theorem~\ref{thm6.1} (and the Markov property of the walk), for that it is enough to prove that for all~$t>0$ and $\BbbP_0$-almost every~$\omega$,
\begin{equation}
\label{6.17a}
\frac{\chi(X_{\lfloor tn\rfloor},\omega)}{\sqrt n}\,\underset{n\to\infty}\longrightarrow\,0
\qquad \text{in $P_{0,\omega}$-probability.}
\end{equation}
Without loss of generality, we need to do this only for~$t=1$. 
By Theorem~\ref{thm6.3a}, the random variable~$X_n$ lies with probability~$1-\epsilon$ in the block $[-M\sqrt n,M\sqrt n]^d\cap\Z^d$, provided~$M$ sufficiently large (with ``large'' depending possibly on~$\omega$).
Using Theorem~\ref{thm6.3a}(1) to estimate $P_{0,\omega}(X_n=x)$ for~$x$ inside this block, we have
\begin{equation}
P_{0,\omega}\bigl(|\chi(X_n,\omega)|>\delta\sqrt n\,\bigr)
\le\epsilon+C(\omega)
\,\frac1{n^{d/2}}\!\!\sum_{\begin{subarray}{c}
x\in\scrC_\infty(\omega)\\|x|\le M\sqrt n
\end{subarray}}
\1_{\{|\chi(x,\omega)|>\delta\sqrt n\}}.
\end{equation}
But Theorem~\ref{thm5.4} tells us that, for all~$\delta,M>0$ and $\BbbP_0$-almost every~$\omega$, the second term tends to zero as~$n\to\infty$. This proves~\eqref{6.17a} and the whole claim.
\end{proofsect}

\subsection{Extension to ``agile'' walk}
\label{sec6.3}\noindent
It remains to prove Theorem~\ref{thm:2ndmainthm} for the ``agile'' version of simple random walk on~$\scrC_\infty$. Since the proof is based entirely on the statement of Theorem~\ref{thm:mainthm}, we will resume a unified treatment of all~$d\ge2$. First we will make the observation that the times of the two walks run proportionally to each other:

\begin{mylemma}
\label{lemma6.3}
Let~$(T_k)_{k\ge0}$ be the stopping times defined in~\eqref{1.6a}. Then for all~$t\ge0$ and $\BbbP_0$-almost every~$\omega$,
\begin{equation}
\frac{T_{\lfloor tn\rfloor}}n\,\underset{n\to\infty}\longrightarrow\,\Theta t,
\qquad P_{0,\omega}\text{\rm-almost surely},
\end{equation}
where
\begin{equation}
\label{E:Theta}
\frac1\Theta=\E_0\Bigl(\frac1{2d}\sum_{e\colon|e|=1}\1_{\{\omega_e=1\}}\Bigr).
\end{equation}
\end{mylemma}

\begin{proofsect}{Proof}
This is an easy consequence of the second part of Theorem~\ref{thm3.1} and the fact that for $\BbbP_0$-almost every~$\omega$ we have~$\tau_x\omega\ne\omega$ once~$x\ne0$. Indeed, let~$f(\omega,\omega')=\1_{\{\omega\ne\omega'\}}$. For~$t=0$ the statement holds trivially so let us assume that~$t>0$. If~$n$ is so large that~$T_{\lfloor nt\rfloor}>0$, we have
\begin{equation}
\frac n{T_{\lfloor tn\rfloor}}=\frac1{T_{\lfloor tn\rfloor}}\sum_{k=1}^{T_{\lfloor tn\rfloor}}f(\tau_{X_{k-1}}\omega,\tau_{X_k}\omega).
\end{equation}
Since~$T_{\lfloor tn\rfloor}\to\infty$ as~$n\to\infty$, by Theorem~\ref{thm3.1} the right hand side converges to the expectation of~$f(\omega,\tau_{X_1}\omega)$ in the annealed measure~$\E_0(P_{0,\omega}(\cdot))$. A direct calculation shows that this expectation equals~$\Theta$.
\end{proofsect}

\begin{proofsect}{Proof of Theorem~\ref{thm:2ndmainthm}}
The proof is based on a standard approximation argument for stochastic processes. Let~$\widetilde B_n(t)$ be as in Theorem~\ref{thm:mainthm} and recall that $\widetilde B_n'(t)$ is a linear interpolation of the values $\widetilde B_n(T_k/n)$ for~$k=0,\dots,n$. The path-continuity of the processes~$\widetilde B_n(t)$ as well as the limiting Brownian motion implies that for every~$\epsilon>0$ there is a~$\delta>0$ such that
\begin{equation}
P_{0,\omega}\Bigl(\,\sup_{\begin{subarray}{c}
t,t'\le T\\|t-t'|<\delta
\end{subarray}}
\bigl|\widetilde B_n(t)-\widetilde B_n(t')\bigr|<\epsilon\Bigl)>1-\epsilon
\end{equation}
once~$n$ is sufficiently large.
Similarly, Lemma~\ref{lemma6.3}, the continuity of~$t\mapsto\Theta t$ and the monotonicity of~$k\mapsto T_k$ imply that for~$n$ sufficiently large,
\begin{equation}
P_{0,\omega}\biggl(\,\sup_{t\le T}\Bigl|\frac{T_{\lfloor tn\rfloor}}n-\Theta t\Bigr|<\delta\biggr)>1-\epsilon.
\end{equation}
On the intersection of these events, the equality~$\widetilde B_n'(k/n)=\widetilde B_n(T_k/n)$ yields
\begin{equation}
\max_{0\le k\le \lfloor Tn\rfloor}\,\bigl|\widetilde B_n'(\ffrac kn)-\widetilde B_n(\ffrac{\Theta k}n)\bigr|<
\epsilon.
\end{equation}
In light of piece-wise linearity this shows that, with probability at least~$1-2\epsilon$, the paths~$t\mapsto\widetilde B_n'(t)$ and~$t\mapsto\widetilde B_n(\Theta t)$ are within a multiple of
$\epsilon$ in the supremum norm of each other. In particular, if~$B_t$ denotes the weak limit of the process~$(B_n(t)\colon t\le T)$, then $(\widetilde B_n'(t)\colon t\le T)$ converges in law to~$(B_{\Theta t}\colon t\le T)$. The latter is an isotropic Brownian motion with diffusion constant~$D'=D\Theta^2$.
\end{proofsect}

\appendix

\section{Heat-kernel upper bounds}
\label{appA}\noindent
Let $(Z_t)_{t\ge0}$ denote the continuous-time random walk which attempts a jump to one of its near\-est-neighbors at rate one (regardless of the number of accessible neighbors). Let~$q_t^\omega(x,y)$ denote the probability that~$Z_t$ started at~$x$ is at~$y$ at time~$t$. In his paper~\cite{Barlow}, Barlow proved the following statement: There exist constants~$C_1,C_2\in(0,\infty)$  and, for each~$x\in\Z^d$, a random variable~$S(x)=S(x,\omega)\in(0,\infty)$ such that for all~$x,y\in\scrC(\omega)$ and all~$t>S(x)$, 
\begin{equation}
\label{A.10}
q_t^\omega(x,y)\le C_1 t^{-d/2}\exp\bigl\{-C_2|x-y|^2/t\bigr\}.
\end{equation}
Moreover,~$S(x)$ has uniformly stretched-exponential tails, i.e.,
\begin{equation}
\label{Sxtail}
\BbbP_0\bigl(S(x)>R\bigr)\le\texte^{-C_3R^{\theta'}},\qquad R\ge1.
\end{equation}
Barlow provides also a corresponding, and significantly harder-to-prove lower bound which requires the additional condition~$t>|x-y|$. However, for~\eqref{A.10}, this condition is redundant.

In the remarks after his Theorem~1, Barlow mentions that appropriate modifications to his arguments yield the corresponding discrete time estimates. Here we present the details of these modifications
which are needed to make our proof of the invariance principles in Theorems~\ref{thm:mainthm} and~\ref{thm:2ndmainthm} complete. Notice that we do not re-prove Barlow's bounds in their full generality, just the absolute minimum necessary for our purposes.

\subsection{Uniform bound}
There will be two kinds of bounds on the heat-kernel as a function of the terminal position of the walk after~$n$ steps: a uniform bound by a constant times~$n^{-d/2}$ and a non-uniform, Gaussian bound on the tails. We begin with the statement of the uniform upper bound:

\begin{myproposition}
\label{propA.1}
Let $d\ge2$ and let $p>p_\cc(d)$. There exists a random variable~$C=C(\omega)$ with
$\BbbP(C<\infty)=1$ such that for all $\omega\in\Omega_0$ and all $x\in\scrC_\infty$,
\begin{equation}
P_{0,\omega}(X_n=x)\le \frac{C(\omega)}{n^{d/2}},
\qquad n\ge1.
\end{equation}
\end{myproposition}

The proof will invoke the isoperimetric bound from Barlow~\cite{Barlow}:

\begin{mylemma}
\label{lemmaA.2}
There exists a constant~$c\in(0,\infty)$ such that for $\BbbP_0$-almost every~$\omega$ and all~$R$ sufficiently large,
\begin{equation}
\frac{|\partial\Lambda|}{|\Lambda|}\ge c|\Lambda|^{-\ffrac1d}
\end{equation}
for all~$\Lambda\subset\scrC_\infty\cap[-R,R]^d$ such that~$|\Lambda|>R^{0.01}$.
\end{mylemma}

\begin{proofsect}{Proof}
This is a consequence of Proposition~2.11 on page~3042, and Lemma~2.13 on page~3045
of Barlow's paper~\cite{Barlow}.
\end{proofsect}

This isoperimetric bound will be combined with the technique of \emph{evolving sets}, developed by Morris and Peres~\cite{Morris-Peres}, whose salient features we will now recall.
Consider a Markov chain on a countable state-space~$V$, let~$p(x,y)$ be the transition kernel and let~$\pi$ be a stationary measure. Let~$Q(x,y)=\pi(x)p(x,y)$ and for each~$S_1,S_2\subset V$, let $Q(S_1,S_2)=\sum_{x\in S_1}\sum_{y\in S_2}Q(x,y)$.
For each set~$S\subset V$ with finite non-zero total measure~$\pi(S)$ we define the \emph{conductance}~$\Phi_S$ by
\begin{equation}
\Phi_S=\frac{Q(S,S^\cc)}{\pi(S)}.
\end{equation}
For sufficiently large $r$, we also define the function
\begin{equation}
\Phi(r)=\inf\bigl\{\Phi_S\colon \pi(S)\le r\bigr\}.
\end{equation}
The following is the content of Theorem~2 in Morris and Peres~\cite{Morris-Peres}: Suppose that~$p(x,x)\ge\gamma$ for some~$\gamma\in(0,\ffrac12]$ and all~$x\in V$. Let~$\epsilon>0$ and~$x,y\in V$. If~$n$ is so large that
\begin{equation}
\label{A.5}
n\ge 1+\frac{(1-\gamma)^2}{\gamma^2}\int_{4[\pi(x)\wedge\pi(y)]}^{4/\epsilon}\frac4{u\Phi(u)^2}\textd u,
\end{equation}
then
\begin{equation}
\label{A.6}
p^n(x,y)\le\epsilon\pi(y).
\end{equation}
Equipped with this powerful result, we are now ready to complete the proof of Proposition~\ref{propA.1}:

\begin{proofsect}{Proof of Proposition~\ref{propA.1}}
First we will prove the desired bound for even times. Fix~$\omega\in\Omega$ and let~$Y_n=X_{2n}$ be the random walk on~$\scrC_\infty(\omega)$ observed only at even times. For each~$x,y\in\scrC_\infty(\omega)$, let us use~$p(x,y)$ to denote the transition probability $P_{x,\omega}(Y_1=y)$. Let~$\pi(x)$ denote the degree of~$x$ on~$\scrC_\infty(\omega)$. Then~$\pi$ is an invariant measure of this chain. Moreover, by our restriction to even times we have~$p(x,x)\ge(2d)^{-2}>0$ and so \twoeqref{A.5}{A.6} can be applied.

By Lemma~\ref{lemmaA.2} we have that~$\Phi_S\ge c\pi(S)^{-\frac1d}$ for some~$c>0$
and all sets~$S$ of the form~$S=\scrC_\infty\cap[-R,R]^d$ for~$R\gg1$. Hence~$\Phi(r)\le c'\, r^{-\frac1d}$ for some finite $c'=c'(\omega)$. Plugging into the integral~\eqref{A.5} and using that~$\pi$ is bounded, we find that if~$n\ge \tilde c\,\epsilon^{-\frac2d}$, then~\eqref{A.6} holds. Here~$\tilde c$ is a positive constant that may depend on~$\omega$. Choosing the minimal~$n$ possible, and applying $p^n(x,y)=P_{x,\omega}(Y_n=y)$, the bound~\eqref{A.6} proves the desired claim for all even times. To extend the result to odd times, we apply the Markov property at time one.
\end{proofsect}

\subsection{Gaussian tails}
Next we will attend to the Gaussian-tail bound. Given the random variables~$S(x,\omega)$ from \twoeqref{A.10}{Sxtail}, define random variables~$N_x=N_x(\omega)$ by
\begin{equation}
N_x = S(x)\vee\sup_{y\colon y\not=x}\frac{S(y)^2}{|y-x|}.
\end{equation}
Here is a restatement of the corresponding bound from Theorem~\ref{thm6.3a}:

\begin{myproposition}
\label{propA.3}
Let $d\ge2$ and $p>p_\cc(d)$. There exist constants~$\c_1,\c_2\in(0,\infty)$ such that for all~$\omega\in\Omega_0$, all $x\in\scrC_\infty(\omega)$, all~$R\ge1$ and all $n>N_x(\omega)$,
\begin{equation}
\sum_{y\colon|y-x|>R}P_{x,\omega}(X_n=y)<\c_1\exp\bigl\{-\c_2R^2/n\bigr\}.
\end{equation}
\end{myproposition}

\begin{proofsect}{Proof}
The proof is an adaptation of Barlow's Theorem~1 to the discrete setting.
Let $(X_n)$ be the discrete time random walk, and let $(Z_t)_{t\ge0}$ be the
continuous time random walk with jumps occurring at rate~$1$, both started at~$x$. We consider the coupling of the two walks such that they make the same moves. We will use~$P$ and~$E$ to denote the coupling measure and the corresponding expectation, respectively.

Let~$n\ge N_x$ and let~$A_n$ be the event that $|X_n-x|>R$. Pick $K>1$ and let
\begin{equation}
I_n=\int_n^{4n}\1_{\{|Z_t-x|>R/K\}}\,\textd t
\end{equation}
be the amount of time in $[n,4n]$ that the walk~$(Z_t)$ spends at distance larger than~$\ffrac RK$ from~$x$. By
the inequality
\begin{equation}
P(A_n)\leq\frac{E(I_n)}{E(I_n|A_n)},
\end{equation}
it suffices to derive an appropriate upper bound on~$E(I_n)$ and a matching lower bound on~$E(I_n|A_n)$. Note that we may assume that~$R\le n$ because otherwise we have~$P(A_n)=0$ and there is nothing to prove.

To derive an upper bound on~$E(I_n)$, we note that for~$t>n$, our choice~$n\ge N_x$ implies~$t>S(x)$. The expectation can then be bounded using \eqref{A.10}:
\begin{equation}
\begin{aligned}
E(I_n) 
&=
\int_n^{4n}\sum_{y\colon|y-x|>R/K}q_t(x,y)\,\textd t 
\\
&\le C_1\int_n^{4n}t^{-d/2}\sum_{x\colon|x|>R/K}\texte^{-C_2|x|^2/t}\,\textd t
\le C_4n\texte^{-C_5\,\frac{R^2}n},
\end{aligned}
\end{equation}
where~$C_4$ and~$C_5$ are constants (possibly depending on~$K$).

It thus remains to prove that, for some constant~$C_6>0$,
\begin{equation}
E(I_n|A_n) \ge C_6\,n.
\end{equation}
To derive this inequality, let us recall that the transitions of~$Z_t$ happen at rate one, and they are independent of the path of the walk. Hence, if~$B_n$ is the event that~$Z_t$ attempted at least~$n$ jumps by time~$2n$, then $P(B_n|A_n)=P(B_n)$ is bounded away from zero for all~$n\ge1$. Therefore, it suffices to prove that~$E(I_n|A_n\cap B_n)\ge C_6 n$. 

Let~$T$ be the first time when the walk~$(Z_t)$ is farther from~$x$ than~$R$. On~$A_n\cap B_n$, this happens before time~$2n$, i.e.,~$T\le 2n$. Let~$Q_R=[-R,R]^d\cap\Z^d$ and $Q_{R/K}=[-\ffrac RK,\ffrac RK]^d\cap\Z^d$. Then for values~$z$ on the external boundary of~$Q_R$---which are those that~$Z_T$ can take---the bound~\eqref{A.10} tells us
\begin{equation}
\label{eq:mir}
\sum_{y\in Q_{R/K}}q_t^\omega(z,y)
\leq C_1\left(\frac{2R}{K}\right)^d \max_{s>0}
\left\{s^{-d/2}\exp(-\tfrac14C_2R^2/s)\right\}\leq C_7K^{-d},
\end{equation}
provided that~$t>S(z)$. But our assumptions~$n\ge N_x$ and $R\le n$ imply~$n\ge S(z)$, and so in light of the fact that~$T\le2n$ on~$A_n\cap B_n$, \eqref{eq:mir} actually holds for all~$t$ such that~$T+t\in[3n,4n]$. 
Plugging~$Z_T$ for~$z$ on the left-hand side and taking expectation gets us an upper bound on $P(Z_t\in Q_{R/K}|A_n\cap B_n)$---with~$t$ now playing the role of~$T+t$. Hence,
\begin{equation}
E(I_n|A_n\cap B_n)\ge\int_{3n}^{4n}P(Z_t\not\in Q_{R/K}|A_n\cap B_n)\ge n(1-C_7K^{-d}).
\end{equation}
Choosing~$K$ sufficiently large, the right-hand side grows linearly in~$n$.
\end{proofsect}

\begin{proofsect}{Proof of Theorem~\ref{thm6.3a}}
Part~(1) is a direct consequence of Proposition~\ref{propA.1}, while part~(2) follows from Proposition~\ref{propA.3} and the fact that if the~$S(x)$ have stretched exponential tails (uniformly in~$x$), then so do the~$N_x$'s.
\end{proofsect}

\section{Some questions and conjectures}
\label{appB}\noindent
While our control of the corrector in~$d\ge3$ is sufficient to push the proof of the functional CLT through, it is not sufficient to provide the \emph{conceptually correct} proof of the kind we have constructed for~$d=2$. However, we do not see any reason why~$d\ge3$ should be different from~$d=2$, so our first conjecture is:

\begin{myconjecture}
\label{conj:cor_small}
Theorem~\ref{thm5.1} is true in all~$d\ge2$.
\end{myconjecture}

Our proof of Theorem~\ref{thm5.1} in~$d=2$ hinged on the fact that the corrector plus the position is a harmonic function on the percolation cluster. Of interest is the question whether harmonicity is an essential ingredient or just mere convenience. Yuval Peres suggested the following generalization of Conjecture \ref{conj:cor_small}:

\begin{myquestion}
\label{ques:yuval}
Let $f\colon\Z^d\to\R$ be a shift invariant, ergodic process on~$\Z^d$ whose gradients are in~$L^1$
and have expectation zero. Is it true that
\begin{equation}
\label{E:B.0}
\lim_{n\to\infty}\frac 1n \max_{x\in \Z^d\cap[-n,n]^d}\,\bigl|f(x)\bigr| = 0
\end{equation}
almost surely?
\end{myquestion}

\noindent
\emph{Update}: The above question, while obviously true in~$d=1$, has a negative answer in all~$d\ge2$. The first counterexample, based on constructions in~\cite{zerner-merkl} and~\cite{BZZ}, was provided to us by Martin Zerner. Later Tom Liggett pointed out the following, embarrassingly simple, counterexample: Let~$f(x)$ be i.i.d.\ with distribution function~$P(f(x)>u)=u^{-d}$ for~$u\ge1$. Then~$(f(x))_{x\in\Z^d}$ is shift-invariant, ergodic, with $f\in L^1$ and the gradients of~$f$ having zero mean, yet $n^{-1}\max_{|x|\le n}|f(x)|$ has a non-trivial distributional limit as~$n\to\infty$.

\smallskip
The harmonic embedding of~$\scrC_\infty$ has been indispensable for our proofs, but it also appears to be a very interesting object in its own right. This motivates many questions about the corrector~$\chi(x,\omega)$. Unfortunately, at the moment it is not even clear what properties make the corrector unique. The following question has been asked by Scott Sheffield:

\begin{myquestion}
Is it true that, for~a.e.~$\omega\in\Omega_0$, there exists only one vector-valued function~$x\mapsto\chi(x,\omega)$ on~$\scrC_\infty(\omega)$ such that $x\mapsto x+\chi(x,\omega)$ is harmonic on~$\scrC_\infty(\omega)$, $\chi(0,\omega)=0$ and~$\chi(x,\omega)/|x|\to0$ as~$|x|\to\infty$?
\end{myquestion}

If this question is answered in the affirmative, we could generate the corrector by its finite-volume approximations (this would also fully justify Fig.~\ref{fig1}). If we restrict ourselves to functions that have the shift-invariance property \eqref{2.12d}, uniqueness can presumably be shown using the ``electrostatic methods'' from, e.g.,~\cite{Golden-Papanicolau}. However, it is not clear whether \eqref{2.12d} holds for the corrector defined by the thermodynamic limit from finite boxes.

As to the more detailed properties of the corrector, for the purposes of the present work one would like to know how $\chi(x,\omega)$ scales with~$x$ and whether it has a well-defined scaling limit. We believe that, in sufficiently high dimension, the corrector is actually tight:

\begin{myconjecture}
\label{conj4}
Let~$d\gg1$. Then for each~$\epsilon>0$ there exists~$K<\infty$ such that $\BbbP_0(\,|\chi(x,\omega)|>K\mid x\in\scrC_\infty)<\epsilon$ for all $x\in\Z^d$.
\end{myconjecture}

It appears that one might be able to prove Conjecture~\ref{conj4} by using Barlow's heat-kernel estimates. To capture the behavior in low dimensions, we make a somewhat wilder guess:

\begin{myconjecture}
\label{conj5}
Let $d\ge1$. Then the law of~$x\mapsto\epsilon^{\frac{2-d}2}\chi(\lfloor\ffrac x\epsilon\rfloor)$ on compact subsets of~$\R^d$ converges weakly (as $\epsilon\downarrow0$) to Gaussian Free Field, i.e., a multivariate Gaussian field with covariance proportional to~$\Delta^{-1}\1$, where~$\Delta$ is the Dirichlet Laplacian on~$\R^d$ and~$\1$ is the $d$-dimensional unit matrix. 
\end{myconjecture}

Here is a heuristic reasoning that led us to these conjectures: Consider the problem of random conductances to avoid problems with conditioning on containment in the infinite cluster. To show the above convergence, we need that for any smooth~$f\colon\R^d\to\R$ with compact support,
\begin{equation}
\label{E:B.2}
\epsilon^d\sum_{x\in\Z^d}(\Delta f)(x\epsilon)\,\epsilon^{\frac{2-d}2}\!\chi(x)
\,\,\overset{\DD}{\underset{\epsilon\downarrow0}\longrightarrow}\,\,
\NN\bigl(0,\sigma^2\Vert\nabla f\Vert_2^2\1\bigr),
\end{equation}
where~$\nabla$ and~$\Delta$ denote the (continuous) gradient and Laplacian, respectively, and where~$\NN(0,C)$ is a mean-zero, covariance-$C$ multivariate normal random variable.
Next we note that the corrector is defined, more or less, as the solution to the equation $\Delta_\textd\chi=-V$, where~$V$ is the local drift and~$\Delta_\textd$ is the relevant generator, which is basically a discrete Laplacian on~$\Z^d$. Thus, if~$g\colon\R^d\to\R$ is smooth with compact support and~$g_\epsilon(x)=g(x\epsilon)$, then
\begin{equation}
\label{E:B.3}
\begin{aligned}
\epsilon^{\frac{d+2}2}\sum_{x\in\Z^d}\epsilon^{-2}(\Delta_\textd & g_\epsilon)(x)\chi(x)
=-\epsilon^{\frac{d-2}2}\sum_{x\in\Z^d}g_\epsilon(x)V(x)
\\
&=\epsilon^{d/2}\sum_{e\colon|e|=1}\sum_{x\in\Z^d}
\frac{g_\epsilon(x+e)-g_\epsilon(x)}\epsilon\,e\,\omega_{(x,x+e)}
\\
&\!\!\!\overset{\DD}{\underset{\epsilon\downarrow0}\longrightarrow}\,\,
\NN\bigl(0,\Vert\nabla g\Vert_2^2\1\bigr)
\end{aligned}
\end{equation}
The convergence statement \eqref{E:B.2} would then follow from \eqref{E:B.3} provided we can replace the ``discretized'' Laplacian $\epsilon^{-2}\Delta_\textd g_\epsilon$ by its continuous counterpart $\Delta g$.

Note that for~$d=1$ and conductances bounded away from zero, Conjecture~\ref{conj5} is actually a theorem.  Indeed, the corrector is a random walk with increments given by reciprocal conductances and so the convergence follows by the invariance principle for random walks. Conjecture~\ref{conj5} suggests that Conjecture~\ref{conj4} applies for~$d\ge3$.

\smallskip
Despite the emphasis on the harmonic embedding of~$\scrC_\infty$, our proofs used, quite significantly, the underlying group structure of~$\Z^d$; e.g., in Sect.~\ref{sec:sec4}. Presumably this will not prevent application of our method to other regular lattices, but for more irregular graphs, e.g., Voronoi percolation in~$\R^d$, significant changes may be necessary. A similar discussion applies to various natural subdomains of~$\Z^d$; for instance, it is not clear how to adapt our proof to random walk on the infinite percolation cluster in the half-space $\{x\in\Z^d\colon x_d\ge0\}$. 

A different direction of generalizations are the models of \emph{long-range} percolation with power-law decay of bond probabilities. Here we conjecture:

\begin{myconjecture}
\label{conj:stable}
Let~$d\ge1$ and consider long-range percolation obtained by adding to~$\Z^d$ a bond between every two distinct sites~$x,y\in\Z^d$ independently with probability proportional to~$|x-y|^{-(d+\alpha)}$. If~$\alpha\in(0,2)$, then the corresponding random walk scales to a symmetric~$\alpha$-stable Levy process in~$\R^d$.
\end{myconjecture}

Note that, according to this conjecture, in~$d=1$, the interval~$\alpha\in(0,2)$ of ``interesting'' exponents is larger than the interval for which an infinite connected component may occur even without the ``help'' of nearest neighbor connections. On the other hand, in dimensions~$d\ge3$, the interval conjectured for stable convergence is strictly smaller than that of ``genuine'' long-range percolation behavior, as defined, e.g., in terms of the scaling of graph distance with Euclidean distance; cf~\cite{Benjamini-Berger,Berger,Biskup}.

\section*{Acknowledgments}
\noindent 
The research of M.B.~was supported by the NSF grant~DMS-0306167. Part of the research was performed while
N.B.~visited ETH-FIM and M.B. visited Microsoft Research in Redmond. We wish to thank these institutions for their hospitality and financial support.
We are also grateful to G.Y.~Amir, A. De Masi, A.~Dembo, P.~Ferrari, T.~Liggett, S.~Olla, Y.~Peres, O.~Schramm, S.~Sheffield, V.~Sidoravicius, A.-S.~Sznitman and M.P.W.~Zerner for interesting and useful discussions at various stages of this project.

\end{document}